\documentclass[11pt]{article}
\usepackage[a4paper,margin=1in]{geometry}
\usepackage{amsmath,amssymb,amsthm,mathtools,mathrsfs}
\usepackage[T1]{fontenc}
\usepackage{lmodern}
\usepackage{microtype}
\usepackage{enumitem}
\usepackage{graphicx} 
\usepackage{soul}
\usepackage{xcolor}
\sethlcolor{yellow}

\usepackage[numbers,sort&compress]{natbib}
\allowdisplaybreaks

\newtheorem{lemma}{Lemma}[section]
\newtheorem{theorem}[lemma]{Theorem}
\newtheorem{definition}[lemma]{Definition}
\newtheorem{corollary}[lemma]{Corollary}
\newtheorem{proposition}[lemma]{Proposition}
\newtheorem{remark}[lemma]{Remark}

\newcommand{\R}{\mathbb R}
\newcommand{\C}{\mathbb C}

\newcommand{\D}{\mathcal D}

\newcommand{\norm}[1]{\left\lVert #1\right\rVert}

\newcommand{\supp}{\operatorname{supp}}

\newcommand{\FL}[2]{\mathcal F L^{#1}_{#2}}
\newcommand{\br}[1]{\langle #1\rangle}
\newcommand{\Sclass}{\mathscr S}
\newcommand{\J}{\mathcal J}
\DeclareMathOperator{\Repart}{Re}
\graphicspath{{results_exp1/}{results_exp2/}{results_exp3/}{results_exp4/}}

\title{Weighted STSAFT-Modulation Spaces and Approximation by Phase-Adapted Neural Networks}

\author{
 Ming Hu$^{1}$ \quad Bei Liu$^{2,*}$\quad
 Bolin Liu$^{1}$ \quad Jie Shen$^{1}$\\[0.6em]
 \parbox{0.92\textwidth}{\centering
   \small $^1$School of Mathematical Sciences and LPMC, Nankai University, Tianjin 300071, China\\
   \small $^2$School of Mathematical Sciences and Institute of Signal Processing and Machine Learning, Tianjin University of Technology, Tianjin 300384, China
 }\\[0.3em]
 \small $^*$Corresponding author: \texttt{liubei1@mail.nankai.edu.cn}\\[0.5em]
 \parbox{0.92\textwidth}{\centering
   \small\textbf{Funding}\\[-0.1em]
   \small This work was supported in part by the National Natural Science Foundation of China (Grant No. 12471131) and the Tianjin Natural Science Foundation (Grant No. 24JCYBJC00480).
 }
}

\date{}

\begin{document}

\maketitle
\begin{abstract}
We develop a transform-adapted approximation framework for shallow neural networks based on the short-time special affine Fourier transform (STSAFT). Replacing the short-time Fourier transform (STFT) in the definition of classical modulation spaces, we work with weighted STSAFT-modulation spaces, building on the existing SAFT-associated modulation-space framework, and use their chirp correspondence with classical weighted modulation spaces through the chirp factorization of the SAFT kernel. We construct an associated neural dictionary and establish an $N$-term Sobolev approximation bound of order $N^{-1/2}$ on bounded one-dimensional domains. A global approximation theorem is also obtained when the spatial translation parameter ranges over the whole real line; a compact-support argument explains why an unrestricted global conclusion cannot in general be inferred from a bounded translation range. For a fixed-scale chirp family $f_\kappa(x)=e^{i\pi\kappa x^2}u(x)$, we combine the global estimate with a phase analysis to obtain the stronger phase-dependent upper-bound factor
$\langle\lambda\rangle^n\langle\kappa-\lambda\rangle^{n+3/2-1/r}.$
For each fixed nonzero Schwartz envelope, matching the dictionary phase reduces the growth of this specified upper bound from order $|\kappa|^{n+3/2-1/r}$ to order $|\kappa|^n$, while preserving the $N^{-1/2}$ rate. The associated mismatch scale is stated as a comparison of upper bounds, not as a threshold for the optimal approximation error. Four numerical experiments examine model comparison, error decay with dictionary size, phase mismatch, and the obstruction caused by restricted translations. We discuss the effective number of fitted output coefficients, including the rank loss of the sine branch at zero phase, and distinguish the numerical diagnostics from the theoretical best-approximation estimates.
\end{abstract}


\section{Introduction}

The Fourier transform provides one of the fundamental tools for analyzing the spectral content of a signal, but its global nature makes it insufficient for describing signals whose frequency content varies with time. This limitation has motivated the development of time--frequency methods, among which the short-time Fourier transform (STFT) plays a central role. By localizing a signal with a window function before applying the Fourier transform, the STFT provides a joint representation in the time--frequency plane and has become a standard tool in harmonic analysis, signal processing, and the study of time-dependent phenomena \cite{Grochenig2001,CorderoRodino2020,Skrettingland2022}. In this framework, localized plane waves describe ordinary oscillatory behavior, while the associated phase-space geometry is determined by translations and modulations. More general signal geometries require more flexible integral transforms. The special affine Fourier transform (SAFT) provides such a generalization by incorporating affine transformations, quadratic-phase factors, and frequency or spatial offsets into a unified framework \cite{T31,T32}. Various special cases of the SAFT include well-known Fourier-type transforms, making it suitable for signals and optical systems involving nontrivial affine phase-space structures \cite{T31}. To introduce temporal localization into the SAFT framework, windowed and short-time variants have been developed. In particular, the windowed SAFT was introduced together with orthogonality relations, inversion formulas, range characterizations, and reconstruction results \cite{T33}, while local and global uncertainty principles for the STSAFT were established in \cite{T34}. Further analytical foundations, including convolution structures, Parseval-type identities, inversion formulas, range characterizations, and Heisenberg-type uncertainty principles, have since been developed \cite{T35}, alongside results on numerical inversion of windowed SAFTs \cite{T36}. These developments indicate that the STSAFT provides a natural phase-space framework for signals whose local oscillations are better described by affine and quadratic-phase structures than by ordinary modulated translations.

Modulation spaces provide a particularly effective functional framework for measuring the joint localization of a function in the time and frequency variables. Originally introduced by Feichtinger and subsequently developed extensively by Feichtinger and Gr\"ochenig, modulation spaces are defined through mixed-norm estimates of the STFT and have become a standard tool in time--frequency analysis \cite{Feichtinger1983,Grochenig2001,GalperinSamarah2004,CorderoRodino2020}. Their importance is reflected in their connections with Gabor expansions, atomic decompositions, pseudodifferential operators, Fourier multipliers, and a wide range of classical function spaces \cite{GalperinSamarah2004,BenyiGrafakosGrochenigOkoudjou2005,FeichtingerNarimani2006}. In parallel, approximation theory for neural networks has developed from the classical universal approximation theorem \cite{T01} toward quantitative theories in which approximation rates depend explicitly on the regularity, complexity, and dimension of the target function class. Barron's seminal work showed that certain Fourier-analytic classes admit dimension-independent approximation rates for shallow networks, while subsequent studies established quantitative results for weighted Fourier--Lebesgue spaces, spectral Barron-type classes, ReLU and ReLU-type networks, and more general activation functions \cite{T02,T03,T04,T05,T07,T08,T19,T23}. The connection between neural networks and nonlinear approximation has further been clarified through variation spaces and dictionary-based formulations \cite{T06,T20,T21,T22,T26}, which interpret shallow networks as sparse nonlinear approximation schemes. Related work has investigated deep-network approximation, mixed smoothness, sampling complexity, optimal rates, and the dependence of network complexity on dimension and regularity \cite{T15,T24,T25,T27,T28,T29,T30}.

Motivated by these developments, we investigate how the geometry of nonstationary signals can be incorporated directly into both the function space and the neural-network architecture. Classical modulation-space approaches use STFT coefficients as the underlying phase-space representation and lead naturally to dictionaries of localized and modulated atoms. Such representations are particularly well suited to signals whose local oscillatory behavior is dominated by translations and ordinary frequency modulation, and are closely related to the interpretation of shallow networks as sparse approximation schemes in variation spaces \cite{T06,T20,T21,T26}. The directly related STFT-based neural approximation framework was developed in \cite{AbdeljawadCordero2026}, while metaplectic neural dictionaries were studied in \cite{AbdeljawadCarioniCordero2026}. More generally, quantitative approximation theory shows that the choice of function class and representation can affect approximation bounds and their dependence on dimension \cite{T02,T03,T19,T24,T27,T28,T29}. These observations motivate a dictionary whose phase-space geometry is adapted to affine and quadratic-phase oscillations.

Specifically, we replace the STFT in the definition of classical modulation spaces with the STSAFT and introduce weighted STSAFT-modulation spaces $M^{p,q}_{m,A}(\mathbb{R})$. Their quasi-norms are determined by weighted mixed-norm estimates of STSAFT coefficients. This construction retains the phase-space viewpoint of modulation spaces \cite{Feichtinger1983,Grochenig2001,CorderoRodino2020} while incorporating the affine and quadratic-phase geometry of the SAFT, in line with the broader role of phase-space geometry and metaplectic transformations in time--frequency analysis \cite{CorderoGiacchi2024}. It also leads to an STSAFT-adapted neural dictionary: by combining the STSAFT inversion formula with a localized activation representation, we obtain atoms that incorporate spatial localization, frequency modulation, and the quadratic phase prescribed by the SAFT parameters. The construction follows the variation-space and dictionary-based approach to shallow-network approximation \cite{T06,T20,T21,T26}.

Our analysis establishes a quantitative Sobolev approximation result on bounded domains. Using the STSAFT inversion formula, a suitable normalization of the dictionary, and variation-space estimates \cite{T06,T20,T26}, we show that, under appropriate assumptions, functions in the proposed spaces admit $N$-term approximation with error of order $N^{-1/2}$ in $W^{n,r}(\Omega)$ for bounded $\Omega\subset\mathbb{R}$. This one-dimensional estimate extends the dictionary-based approximation framework \cite{T06,T20,T21} to a representation adapted to STSAFT geometry and is consistent with established quantitative results for shallow networks \cite{T02,T03,T19,T20,T21}. Appropriate embeddings further yield related bounds for weighted Fourier--Lebesgue and A-Barron classes, connecting the present approach with existing Fourier-analytic and weighted-space approximation theories \cite{T02,T03,T04,T05,T07,T08,T19,T23,T24}.

We also quantify the dependence of the approximation bound on the quadratic phase encoded by the dictionary. For a fixed-scale SAFT family and chirp targets $f_\kappa(x)=e^{i\pi\kappa x^2}u(x)$, the phase parameter $\lambda$ determines a trade-off between the Sobolev cost of restoring the dictionary phase and the representation cost associated with the residual phase $\kappa-\lambda$. A stronger input estimate follows by applying the whole-space approximation theorem with unrestricted translations and then restricting the approximants to the bounded interval. With $s_*=n+1-1/r$, this gives the phase-dependent factor $\langle\lambda\rangle^n\langle\kappa-\lambda\rangle^{n+3/2-1/r}$. For each fixed nonzero Schwartz envelope, full phase matching reduces the growth of this specified upper bound from order $|\kappa|^{n+3/2-1/r}$ to order $|\kappa|^n$, without changing the $N^{-1/2}$ rate. The resulting mismatch scale is a comparison of these particular upper bounds; it is not a threshold for the optimal approximation error or for empirical training errors. The analysis is motivated by the broader role of structured function classes and dictionaries in efficient approximation \cite{T02,T03,T19,T20,T21,T24,T27,T28,T29}.

Finally, four numerical experiments examine the practical behavior of the proposed construction. Under approximately matched nominal full-network parameter counts, we compare the STSAFT-inspired dictionary with a shallow ReLU spline on stationary and chirp-type targets, assessing both function approximation and derivative reconstruction. Additional experiments investigate error decay as the dictionary grows, the effect of matching or mismatching the quadratic phase, and the support obstruction that arises when a compactly supported window is combined with a restricted translation range. The results indicate that phase adaptation can improve approximation for the tested chirp targets, particularly in derivative reconstruction. These empirical comparisons are qualified by the effective output dimension, which is smaller at zero phase because the sine branch vanishes.

The contributions are organized as follows. First, we place the approximation problem in the weighted STSAFT-modulation framework and make explicit its chirp correspondence with classical modulation spaces; the underlying SAFT-associated modulation spaces themselves have been studied previously in \cite{T37}. Second, we construct an STSAFT-adapted neural dictionary and prove a local Sobolev approximation theorem with an $N^{-1/2}$ rate, together with conditional Fourier--Lebesgue and A-Barron consequences. Third, we establish a whole-space approximation theorem using all spatial translations and identify why a bounded translation restriction does not in general support the same global reconstruction argument. Fourth, we combine the global estimate with a fixed-envelope chirp asymptotic to obtain a stronger phase-dependent upper bound and its associated mismatch scale. Fifth, we provide numerical diagnostics for model comparison, dictionary-size effects, phase matching, and restricted translations, while distinguishing nominal parameter counts from the effective output dimension. The local STFT-based neural approximation framework of \cite{AbdeljawadCordero2026} and the metaplectic dictionary framework of \cite{AbdeljawadCarioniCordero2026} provide the closest points of comparison; the present results concern the STSAFT-adapted dictionary, the all-translation repair, and the quantitative fixed-envelope phase estimate.

\section{Preliminaries}

\subsection{The Short-Time Fourier Transform}

In the field of signal analysis and time-frequency methods, analyzing a signal exclusively in the time domain or frequency domain is often inadequate. The short-time Fourier transform (STFT) is adopted to characterize the temporal evolution of signal frequency components. While the conventional Fourier transform only provides a global frequency representation of signals, the STFT employs a window function to temporally localize the signal prior to Fourier transformation. This yields a bivariate function that simultaneously describes the time and frequency characteristics of the signal. If we introduce the translation $T_x$ and modulation $M_\omega$ operators, namely
$T_x f(t) = f(t - x), \quad M_\omega f(t) = e^{2\pi i \omega\cdot t} f(t),$
the STFT of a signal $f\in L^2(\mathbb R^d)$ with respect to a non-zero window $g\in L^2(\mathbb R^d)$ is given by
\[
(V_g f)(x,\omega) = \langle f, M_\omega T_x g\rangle_{L^2} = \mathcal{F}(f\cdot T_x \overline{g})(\omega)
= \int_{\mathbb R^d} f(y)\,\overline{g(y - x)} \,e^{-2\pi i y\cdot \omega} dy. \tag{2.1} \label{2.2}
\]

Fix $\gamma \in L^2(\mathbb R^d)$. The STFT analysis operator $V_\gamma$ is a bounded map from $L^2(\mathbb R^d)$ into (not onto) $L^2(\mathbb R^{2d})$. Its adjoint $V_\gamma^*$ is given by
$V_\gamma^* F = \int_{\mathbb R^{2d}} F(x,\omega) M_\omega T_x \gamma \,dx d\omega.$
The operator $V_\gamma^*$ is a bounded operator from $L^2(\mathbb R^{2d})$ onto $L^2(\mathbb R^d)$. For $F = V_g f$, with $g,\gamma \in L^2(\mathbb R^d)$, $(g,\gamma) \neq 0$, the inversion formula is given by
\[f = \frac{1}{(g,\gamma)} V_\gamma^* V_g f,\quad f \in L^2(\mathbb R^d).\]

\subsection{The Short-Time Special Affine Fourier Transform}

The Special Affine Fourier Transformation (SAFT), first introduced in \cite{T31}, constitutes an integral transformation in phase space. This transform is associated with general inhomogeneous lossless linear mappings on phase space and relies on six parameters that do not depend on phase-space coordinates. It maps the position $x$ and the wave number $k$ into
\begin{equation}
\begin{bmatrix}
x' \\
k'
\end{bmatrix}
=
\begin{bmatrix}
a & b \\
c & d
\end{bmatrix}
\begin{bmatrix}
x \\
k
\end{bmatrix}
+
\begin{bmatrix}
p \\
q
\end{bmatrix}.\tag{2.2}
\label{eq:saft_phase_map}
\end{equation}
where $\{a,b,c,d,p,q\}$ are considered to be real-valued with
$ad-bc = 1.$

This transformation, which can model many general optical systems, maps any convex body into another convex body and guarantees that the area of the body is preserved by the transformation. The SAFT offers a unified viewpoint of known signal processing transformations as well as optical operations on light waves. 

The definition of SAFT is given as follows:
\begin{equation}
\begin{aligned}
F(w) &= \widehat{f}_{\mathbf{A}}(w) = \int_{\mathbb R} K_{A}(t,w) f(t) dt \quad \text{(SAFT of $f(t)$)}\\
&= \frac{1}{\sqrt{2\pi |b|}} \int_{\mathbb R} \exp\left( \frac{j}{2b} \big(at^2 + dw^2 - 2tw + 2pt + 2(bq-dp)w\big) \right) f(t) dt,
\end{aligned}
\label{eq:saft_def}\tag{2.3}
\end{equation}
where $\mathbf{A}$ stands for the six parameters $(a,b,c,d,p,q)$, and
\[
K_{A}(t,w) = \frac{1}{\sqrt{2\pi |b|}} \exp\left( \frac{j}{2b} \big(at^2 + dw^2 - 2tw + 2pt + 2(bq-dp)w\big) \right).\tag{2.4}\label{2.8}
\]  

The inversion formula for the SAFT is easily shown to be
\begin{equation}
f(t) = \frac{1}{\sqrt{2\pi |b|}} \int_{\mathbb R} F(w) \exp\left( -\frac{j}{2b} \big(dw^2 + at^2 - 2tw + 2w(bq-dp) + 2pt\big) \right) dw.\tag{2.5}
\label{eq:saft_inverse}
\end{equation}

\begin{definition}[STSAFT\cite{T33}]
Let $g\in L^2(\mathbb{R})$ be a window function with
$
\|g\|_{L^2} = 1.
$and

\[
A=
\begin{bmatrix}
a & b & p\\
c & d & q
\end{bmatrix},
\]
be the parameters satisfying $a,b(\neq 0),c,d,p,q\in\mathbb{R}$, and
$ad-bc=1$.

The STSAFT of a complex-valued signal
$f(\tau)\in L^2(\mathbb{R};\mathbb C)$
with respect to $g$ is defined by
\[
S_g^{A}f(x,\omega)
=
\int_{\mathbb{R}}
f(\tau)\overline{g(\tau-x)}
K_{A}(\tau,\omega)\,d\tau ,
\tag{2.6} \label{2.11}
\]
where $K_{A}(\tau,\omega)$ is given by \eqref{2.8}. For tempered distributions, the transform is understood in the standard weak sense (with a Schwartz window), equivalently through the chirp--STFT correspondence established below.
\end{definition}

%
%
%

\subsection{Function Space}

This section compiles definitions and core properties of function spaces used for our analysis. We cover weighted Fourier–Lebesgue spaces, Barron spaces, modulation spaces (including embeddings), and classical Sobolev spaces. They serve as the analytical setting for our approximation results. Most spaces discussed below are formulated using weight functions. To facilitate the discussion, we first present the weight class considered in this section.

\noindent\textbf{Weight Functions.} 
Let $v$ be a continuous, positive, and submultiplicative weight function on
$\mathbb{R}^{d}$, namely, it satisfies
$v(z_1+z_2)\leq v(z_1)v(z_2),z_1,z_2\in\mathbb{R}^{d}.$
 
A function $m$ is said to belong to the class $\mathcal{M}_v(\mathbb{R}^d)$
if it is positive and continuous, and fulfills the following
$v$-moderateness property:
$m(z_1+z_2)\leq C v(z_1)m(z_2),z_1,z_2\in\mathbb{R}^{d},$
where $C>0$ is a constant independent of $z_1$ and $z_2$.

In this work, we mainly consider polynomially growing weights on
$\mathbb{R}^{n}$, where $n=d$ or $n=2d$, which are defined by
\[
v_s(z)=\langle z\rangle^s,
\qquad z\in\mathbb{R}^{n},
\tag{2.7} \label{2.14}
\]
with
$\langle z\rangle=(1+|z|^2)^{1/2}.$

The corresponding tensor product weights on $\mathbb{R}^{2d}$ are given by
\[
(v_s\otimes 1)(x,\xi)=(1+|x|^2)^{s/2},
\qquad
(1\otimes v_s)(x,\xi)=(1+|\xi|^2)^{s/2},
\qquad x,\xi\in\mathbb{R}^{d}.  
\]

It is worth noting that when $s<0$, the weight function $v_s$ is
$v_{|s|}$-moderate.

For two weights $m_1$ and $m_2$ defined on $\mathbb{R}^{d}$, their tensor
product is introduced by
\[
(m_1\otimes m_2)(x,\xi)
=
m_1(x)m_2(\xi),
\qquad x,\xi\in\mathbb{R}^{d}. 
\]

\subsubsection{Weighted Lebesgue And Weighted Fourier-Lebesgue Spaces}

Let $0 < p \leq \infty$ and let $m : \mathbb{R}^d \to (0,\infty)$ be a weight function. The weighted Lebesgue space $L_m^p(\mathbb{R}^d)$ consists of all measurable functions $f : \mathbb{R}^d \to \mathbb{C}$ such that the following (quasi-)norm
\[
\|f\|_{L_m^p(\mathbb{R}^d)} :=
\begin{cases}
\left(\int_{\mathbb{R}^d} |f(x)|^p m(x)^p dx\right)^{1/p}, & 0 < p < \infty,\\[4pt]
\operatorname{ess\,sup}_{x\in\mathbb{R}^d} |f(x)| m(x), & p = \infty,
\end{cases}    
\]
is finite.

Similarly, for $0 < p,q \leq \infty$, and $F : \mathbb{R}^{2d} \to \mathbb{C}$ measurable, we set
\[
\|f\|_{L_m^{p,q}(\mathbb{R}^{2d})} :=
\left(\int_{\mathbb{R}^d} \left(\int_{\mathbb{R}^d} |F(x,y)|^p m(x,y)^p dx\right)^{\frac{q}{p}} dy\right)^{\frac{1}{q}},  
\]
where $m$ is a weight function on $\mathbb{R}^{2d}$.

The weighted Fourier-Lebesgue spaces $\mathcal{F}L_s^p(\mathbb{R}^d)$ are defined in terms of the weighted integrability of the Fourier transform.

\begin{definition}[Weighted Fourier-Lebesgue Spaces]
Let $0 < p \leq \infty$ and $s \in \mathbb{R}$. The weighted Fourier-Lebesgue space $\mathcal{F}L_s^p(\mathbb{R}^d)$ is defined by
\[
\mathcal{F}L_s^p(\mathbb{R}^d) =
\left\{ f \in \mathscr{S}'(\mathbb{R}^d) : \|f\|_{\mathcal{F}L_s^p} := \|v_s \hat{f}\|_{L^p(\mathbb{R}^d)} < \infty \right\},
\label{2.20}
\tag{2.8} 
\]
where $v_s$ is defined in \eqref{2.14}.
\end{definition}

\subsubsection{Barron Spaces}

Barron spaces, originally introduced in the pioneering works of Barron\cite{T19}, and subsequently generalized in various studies\cite{T02,T23,T28}, offer a Fourier-based framework for characterizing functions that can be efficiently approximated by shallow neural networks.

\begin{definition}[Barron Norm and Barron Space]
For $s\in\mathbb{R}$, we define the Barron space as
\[
B_s(\mathbb{R}^d) = \left\{ f\in\mathscr{S}'(\mathbb{R}^d): \|f\|_{B_s}<\infty \right\}, 
\]
where the Barron norm of $f$ is defined as
$\|f\|_{B_s} = \int_{\mathbb{R}^d} (1+|\xi|)^s |\hat{f}(\xi)| \,d\xi.$
\end{definition}

Since $(1+|\xi|)^s\asymp v_s(\xi)$, the spectral Barron space defined above agrees with $\mathcal FL^1_{v_s}$ with equivalent norms; see \cite{T19,T28}.

\subsubsection{Modulation Spaces}

First introduced by Feichtinger in reference \cite{Feichtinger1983} and further extended in subsequent studies such as \cite{T27}, modulation spaces have become a well-established subject within time-frequency analysis. 

Let $g\in\mathcal{S}(\mathbb R^d)$ denote a non-zero window function, $m\in\mathcal{M}_v$, and $0<p,q\le\infty$. The modulation space $M_m^{p,q}(\mathbb R^d)$ consists of all tempered distributions $f\in\mathcal{S}'(\mathbb R^d)$ satisfying
\[
\|f\|_{M_m^{p,q}}=\|V_g f\|_{L_m^{p,q}}=\left(\int_{\mathbb R^d}\left(\int_{\mathbb R^d}|V_g f(x,\omega)|^p m(x,\omega)^p dx\right)^{q/p}d\omega\right)^{1/q}<\infty, 
\]
with standard conventions adopted for the limiting cases $p=\infty$ or $q=\infty$. The STFT $V_g f$ is defined according to equation \eqref{2.2}. We adopt simplified notations: $M_m^p(\mathbb R^d)$ stands for $M_m^{p,p}(\mathbb R^d)$, and $M^{p,q}(\mathbb R^d)$ is used when $m\equiv 1$.

The space $M^{p,q}(\mathbb R^d)$ forms a Banach space for $p,q\ge1$, while it becomes a quasi-Banach space in other scenarios. Up to norm equivalence, its (quasi-)norm is independent of the particular choice of window function $g$, for nonzero Schwartz windows and the polynomial weights used here. Here we use nonzero Schwartz windows and polynomial weights; no unconditional extension to all $M^1$ windows is asserted. 

Modulation spaces satisfy the following inclusion chain: if $0 < p_1 \le p_2 \le \infty$, $0 < q_1 \le q_2 \le \infty$ and $m_1,m_2$ weights in $\mathbb R^{2d}$ which satisfy $m_2 \lesssim m_1$, then
\[
\mathscr{S}(\mathbb R^d) \hookrightarrow M_{m_1}^{p_1,q_1}(\mathbb R^d) \hookrightarrow M_{m_2}^{p_2,q_2}(\mathbb R^d) \hookrightarrow \mathscr{S}'(\mathbb R^d). 
\]

\begin{theorem} (\cite{guo2016}, Theorem 4.11). Let $0 < p_j,q_j \leq \infty$, $s_j,t_j \in \mathbb R$, for $j=1,2$, and consider the polynomial weights $v_{t_j},v_{s_j}$ defined as in \eqref{2.14}. Then
\[
M_{v_{t_1}\otimes v_{s_1}}^{p_1,q_1}(\mathbb R^d) \hookrightarrow M_{v_{t_2}\otimes v_{s_2}}^{p_2,q_2}(\mathbb R^d) 
\]
if the following two conditions hold:

\noindent(i) $(p_1,p_2,t_1,t_2)$ satisfies one of the following:
\begin{align*}
\tag{C1} \frac{1}{p_2} &\leq \frac{1}{p_1},\quad t_2 \leq t_1,\\
\tag{C2} \frac{1}{p_2} &> \frac{1}{p_1},\quad \frac{1}{p_2}+\frac{t_2}{d} < \frac{1}{p_1}+\frac{t_1}{d};
\end{align*}

\noindent(ii) $(q_1,q_2,s_1,s_2)$ satisfies either (C1) or (C2) with $p_j$ replaced by $q_j$ and $t_j$ replaced by $s_j$, respectively.
\end{theorem}

Based on the definition of the STSAFT in \eqref{2.11}, we can define the weighted STSAFT-modulation spaces $M_{v_t \otimes v_s,A}^{p,q}$ as the collection of all functions 
$f\in \mathscr{S}'(\mathbb{R}^d)$ satisfying 
$\|f\|_{M_{v_t \otimes v_s,A}^{p,q}}<\infty,$
where 
$\|f\|_{M_{v_t \otimes v_s,A}^{p,q}} = \|S_g^A f\|_{L_{v_t \otimes v_s}^{p,q}}.$

\begin{lemma}
Let
$0<p_j,q_j\le\infty, t_j,s_j\in\mathbb R, j=1,2.$
Assume that$(p_1,p_2,t_1,t_2)$
satisfies either \emph{(C1)} or \emph{(C2)} of Theorem~2.4, and that
$(q_1,q_2,s_1,s_2)$
satisfies the corresponding analogous condition. Then
\[
M^{p_1,q_1}_{v_{t_1}\otimes v_{s_1},A}(\mathbb R)
\hookrightarrow
M^{p_2,q_2}_{v_{t_2}\otimes v_{s_2},A}(\mathbb R).     
\]
\end{lemma}

\begin{proof}
Define
$ (\mathcal C_Af)(t)=e^{i(at^2+2pt)/(2b)}f(t),\qquad
 \Psi_A(w)=\frac{d}{2b}w^2+\left(q-\frac{dp}{b}\right)w.$
Thus $\Psi_A$ is a real-valued phase. Directly from the kernel,
\begin{equation}
 S_\varphi^Af(y,w)=\frac{e^{i\Psi_A(w)}}{\sqrt{2\pi|b|}}
 V_\varphi(\mathcal C_Af)\left(y,\frac{w}{2\pi b}\right).
 \tag{2.9}\label{eq:STSAFT-STFT-relation}
\end{equation}
The change of variables $w=2\pi b\xi$, with absolute Jacobian $2\pi|b|$, and
$v_s(2\pi b\xi)\asymp v_s(\xi)$ give
\begin{equation}
 \|f\|_{M^{p,q}_{v_t\otimes v_s,A}}
 \asymp\|\mathcal C_Af\|_{M^{p,q}_{v_t\otimes v_s}}.
 \tag{2.10}\label{eq:STSAFT-classical-norm-equivalence}
\end{equation}
This also holds with the usual essential-supremum convention at infinite exponents. Applying the classical weighted embedding in Theorem~2.4 (\cite[Theorem 4.11]{guochen2016}, with dimension one) to $\mathcal C_Af$ proves the claimed STSAFT embedding.
\end{proof}

\subsubsection{Potential Sobolev Spaces}

Let $s\in\mathbb R$ and $1<r<\infty$. The Bessel potential space $W^{s,r}(\mathbb R^d)$ consists of all tempered distributions $f\in\mathscr{S}'(\mathbb R^d)$ obeying
\[
\|f\|_{W^{s,r}}:=\left\|\mathcal{F}^{-1}\big(v_s \hat{f}\big)\right\|_{L^r(\mathbb R^d)}<\infty, 
\]
where $v_s$ denotes the weight introduced in formula \eqref{2.14}. An equivalent characterization reads
\[
W^{s,r}(\mathbb R^d)=\big\{f\in\mathscr{S}'(\mathbb R^d):\langle D\rangle^s f\in L^r(\mathbb R^d)\big\},  
\]
with the Bessel potential operator $\langle D\rangle^s$ defined via
$\langle D\rangle^s f:=\mathcal{F}^{-1}\big(v_s \hat{f}\big).$

For $s=n\in\mathbb Z_+$ and $1<r<\infty$, this space agrees with the weak-derivative Sobolev space with equivalent norms. The Fourier--Lebesgue space $\mathcal FL^r_{v_s}$ is the inverse Fourier image of the weighted space $L^r_{v_s}$, not the Fourier image of a Bessel potential space; see \cite{pilipovic2010}. Moreover, Toft established embedding theorems connecting Sobolev spaces and modulation spaces\cite{toft2004}.

\noindent\textbf{Proposition 5.} Suppose $s\in\mathbb R$, $1\le p,q\le\infty$ and $1<r<\infty$. If $q\le p\le r\le q'$, then
\[
M_{1\otimes v_s}^{p,q}(\mathbb R^d)\hookrightarrow W^{s,r}(\mathbb R^d),   
\]
where the inclusion mapping is continuous.

\subsection{Variation Space and Maurey's Sampling Result}

Let $\mathcal{B}$ denote a complex Banach space (realified when applying the type-2 sampling theorem), and let $\mathcal{D}\subset \mathcal{B}$ stand for a collection of nonzero elements, referred to as a dictionary (namely a set of atoms). In the context of general nonlinear approximation, the ordering of elements within $\mathcal{D}$ carries no significance; only the selection of atoms and their corresponding coefficients matter. For $N\in\mathbb{N}$ and $M>0$, we define the nonlinear manifold of $N$-term approximants with $\ell_1$ regularization associated with the dictionary $\mathcal{D}$ as follows:
\[
\Sigma_{N,M}(\mathcal{D}) := \left\{\sum_{j=1}^N a_j h_j : a_j\in\mathbb C,\ h_j\in\mathcal{D},\ \sum_{j=1}^N |a_j| \le M\right\}. \tag{2.11}
\]
Eliminating the $\ell_1$ regularization constraint yields the $N$-term nonlinear manifold
\[
\Sigma_N(\mathcal{D}) := \bigcup_{M>0}\Sigma_{N,M}(\mathcal{D}).  \tag{2.12}
\]

Accordingly, $\Sigma_{N,M}(\mathcal{D})$ collects all linear combinations formed by at most $N$ atoms selected from $\mathcal{D}$ under the $\ell^1$ regularization framework, while $\Sigma_N(\mathcal{D})$ imposes no regularization. Given a target element in the Banach space $\mathcal{B}$, this nonlinear manifold is employed to construct the optimal atomic combination to approximate the target with minimal error. To formalize this idea, one notices that every bounded linear combination of atoms can be rescaled into a convex combination. This property enables us to quantify the complexity of an element via the minimal scaling factor needed to place it inside the convex hull of the dictionary.

\begin{definition}
Let $\mathcal{B}$ be a Banach space and $\mathcal{D} \subseteq \mathcal{B}$ be a dictionary. Then for $f \in \mathcal{B}$, the variation gauge of $\mathcal{D}$ is defined as
\[
\|f\|_{\mathcal{K}(\mathcal{D})} := \inf\big\{ c > 0 : f/c \in \overline{\mathrm{conv}}(\mathcal{D}^{\mathrm{ph}}) \big\},  \tag{2.13}
\]
where $\mathcal{D}^{\mathrm{ph}}=\{e^{i\theta}h:h\in\mathcal D,\ 0\le\theta<2\pi\}$ and the closed convex hull is taken in the underlying real Banach space. The corresponding variation space is then the set of functions with finite variation norm
\[
\mathcal{K}(\mathcal{D}) := \big\{ f \in \mathcal{B} : \|f\|_{\mathcal{K}(\mathcal{D})} < \infty \big\}.  
\]
\end{definition}

In general, the displayed functional is a gauge (possibly an extended-valued seminorm), not automatically a norm for an arbitrary dictionary. It is a norm on its finite domain when the closed, balanced convex hull of the phase-expanded dictionary is bounded and absorbing in the underlying real Banach space. In the arguments below we use only the gauge and its representing-measure upper bound; no unverified nondegeneracy of the dictionary is required.

This gauge serves as a crucial quantity in determining the approximation rates achieved by sparse representation techniques. 
In \cite{T26}, a generalized version of Maurey’s approximation theorem was developed for functions in the variation space associated with a dictionary, yielding the following approximation result:

\begin{proposition}[Approximation Rate in Type-2 Banach Spaces]
Let $\mathcal{B}$ be a type-2 Banach space and $\mathcal{D} \subset \mathcal{B}$ be a dictionary with $K_{\mathcal{D}} := \sup_{d\in\mathcal{D}} \|d\|_{\mathcal{B}} < \infty$. Then for $f \in \mathcal{K}(\mathcal{D})$, we have
\[
\inf_{f_N\in\Sigma_{N,M_f}(\mathcal{D})} \|f - f_N\|_{\mathcal{B}} \leq 4C_{2,\mathcal{B}} K_{\mathcal{D}} \|f\|_{\mathcal{K}(\mathcal{D})} N^{-\frac12} \tag{2.14}
\]
with $M_f = \|f\|_{\mathcal{K}(\mathcal{D})}$. For complex coefficients, apply the real type-2 result to $\mathcal D^{\mathrm{ph}}$; the unit phases can be absorbed into the coefficients.
\end{proposition}

Related extensions of Maurey's approximation theorem can be found in \cite[Section 8]{T22}. 
We use only the following sufficient condition and upper bound for a given representing measure. The exact measure characterization in \cite[Lemma 3]{T06} requires additional hypotheses, in particular compactness of the dictionary in that result.

\begin{proposition}[\cite{T06}, Lemma 3]
Let $\mathcal{B}$ be a Banach space and suppose that $\mathbb{D} \subset \mathcal{B}$ is compact. Then $f \in \mathcal{K}(\mathcal{D})$ if there exists a Borel measure $\mu$ on $\mathcal{D}$ such that
\[
f = \int_{\mathcal{D}} i_{\mathcal{D}\to\mathcal{B}} \,d\mu.   \tag{2.15}
\]
Moreover,
\[
\|f\|_{\mathcal{K}(\mathcal{D})} = \inf\left\{ \|\mu\| : f = \int_{\mathcal{D}} i_{\mathcal{D}\to\mathcal{B}} \,d\mu \right\},   \tag{2.16}
\]
where the infimum is taken over all Borel measures $\mu$ defined on $\mathcal{D}$, and $\|\mu\|$ is the total variation of $\mu$.
\end{proposition}


\section{Convergence Rate for Approximation of STSAFT-modulation spaces} 
\label{chpt:mainresult}

In this section, we derive several approximation theorems for shallow neural networks targeting functions belonging to weighted modulation spaces $M_m^{p,q}(\mathbb R^d)$, with approximation errors measured in diverse norms.

The activation representation used below is the STFT representation in
\cite{AbdeljawadCordero2026}. For a nonzero real-valued window
$\phi\in\mathscr S(\mathbb R)$ and $V_\phi\sigma(t_0,\tau_0)\ne0$,
$\tau_0\ne0$, a change of variables gives
\[
e^{2\pi i\eta x}=
\frac{1}{V_\phi\sigma(t_0,\tau_0)}
\int_{\mathbb R}\sigma(\eta x/\tau_0+\beta)
\phi(\eta x/\tau_0+\beta-t_0)e^{-2\pi i\beta\tau_0}\,d\beta.
\tag{3.1}
\]
The integral is absolutely convergent since $\sigma$ is bounded and
$\phi$ is integrable.

\subsection{Local Approximation}

\begin{theorem}[Local Approximation in STSAFT-Modulation Spaces]
\label{thm:STSAFT-local-approximation}
Let $n\in\mathbb N$, $0<p<\infty$, $0<q\leq 2\leq r<\infty$,
and let $\Omega\subset\mathbb R$ be a bounded domain. Let
\[
A=
\begin{bmatrix}
a & b & p_0\\
c & d & q_0
\end{bmatrix},
\qquad
b\neq0,
\qquad
ad-bc=1,
\]
be a fixed set of SAFT parameters. Define
\[
\Phi_A(x)
=
\frac{a}{2b}x^2+\frac{p_0}{b}x,
\qquad
\Psi_A(w)
=
\frac{d}{2b}w^2
+
\left(q_0-\frac{dp_0}{b}\right)w, 
\qquad
\eta_A(w)=\frac{w}{2\pi b}. \tag{3.2}
\]

Let $\sigma\in W^{k,\infty}(\mathbb R)\setminus\{0\}, k\ge n,$
and 
\[
\phi\in\mathscr{S}(\mathbb R)\setminus\{0\},
\qquad
g\in\mathscr{S}(\mathbb R)\setminus\{0\},
\qquad
\|g\|_{L^2(\mathbb R)}=1.
\]
As in the activation representation above, take $\phi$ real-valued.
Assume that there exist $t_0\in\mathbb R$ and $\tau_0\neq0$ such that
$V_\phi\sigma(t_0,\tau_0)\neq0.$

For a weight
$m_A(y,w)
=
v_{s_1}(y)v_{s_2}(\eta_A(w)),$
where
\[
s_1
\begin{cases}
=0, & 0<p\le1,\\[1mm]
>1/p', & p>1,
\end{cases}
\qquad
s_2
\begin{cases}
=n+1, & 0<q\le1,\\[1mm]
>n+1+1/q', & q>1,
\end{cases}
\]
consider the STSAFT-modulation space
\[
M_{m_A,A}^{p,q}(\mathbb R)
=
\left\{
f\in\mathscr{S}'(\mathbb R):
\|f\|_{M_{m_A,A}^{p,q}}<\infty
\right\},   \tag{3.3}
\]
with
$\|f\|_{M_{m_A,A}^{p,q}}
=
\|S_g^Af\|_{L_{m_A}^{p,q}(\mathbb R^2)}.  $

Define the SAFT-adapted dictionary $\mathcal{D}_A
=
\left\{
\rho_A(\cdot,y,w,\beta):
(y,w,\beta)\in\mathbb R^3
\right\},$
where
\[
\begin{aligned}
\rho_A(x,y,w,\beta)
:={}&
e^{-i\Phi_A(x)}
\sigma\left(
\frac{\eta_A(w)x}{\tau_0}+\beta
\right)
\\
&\times
\phi\left(
\frac{\eta_A(w)x}{\tau_0}
+\beta-t_0
\right)
g(x-y).  
\end{aligned}
\]

Then, for every
$f\in M_{m_A,A}^{p,q}(\mathbb R)$, there exists a constant
$C>0$, depending on
$A,\Omega,n,r,p,q,\sigma,\phi$ and $g$, but independent of
$f$ and $N$, such that
\[
\inf_{f_N\in\Sigma_N(\mathcal{D}_A)}
\|f-f_N\|_{W^{n,r}(\Omega)}
\le
C
N^{-1/2}
|\Omega|^{1/r}
\|f\|_{M_{m_A,A}^{p,q}}.  \tag{3.4}\label{eq3.10}
\]  \label{thm:main}
\end{theorem}

\begin{proof}
We use the STFT approximation mechanism of \cite{AbdeljawadCordero2026}
and record the changes induced by the SAFT phase.
The STFT inversion formula below is understood in its weak, distributional
form for tempered distributions; no additional assumption $f\in L^2$ is
made. The weighted coefficient estimate and the subsequent Bochner bound
justify its realization in $W^{n,r}(\Omega)$. The kernel factorization gives
\[
S_g^Af(y,w)=\frac{e^{i\Psi_A(w)}}{\sqrt{2\pi|b|}}
 V_g(e^{i\Phi_A}f)(y,\eta_A(w)). 
\]
STFT inversion, $\|g\|_2=1$, and the activation representation yield
\[
f(x)=C_{A,\sigma,\phi}\int_{\mathbb R^3}
 \rho_A(x,y,w,\beta)e^{-i\Psi_A(w)}e^{-2\pi i\beta\tau_0}
 S_g^Af(y,w)\,dy\,dw\,d\beta,\tag{3.5}\label{eq3.14}
\]
where $C_{A,\sigma,\phi}=(\sqrt{2\pi|b|}\,
V_\phi\sigma(t_0,\tau_0))^{-1}$.

Put $R_\Omega=\sup_{x\in\Omega}|x|$, choose $s<-1$, and set
\[
\vartheta_A(w,\beta)=v_n(\eta_A(w))
 v_s\big((|\beta|-R_\Omega|\eta_A(w)|/|\tau_0|)_+\big),
\qquad \widetilde\rho_A=\rho_A/\vartheta_A.
\]
Let $\widetilde{\mathcal D}_A$ be the dictionary of these normalized atoms.
For $r_0=R_\Omega|\eta_A(w)|/|\tau_0|$,
\[
\int_{\mathbb R}\vartheta_A(w,\beta)\,d\beta
=2v_n(\eta_A(w))\left(r_0+\int_0^\infty(1+t^2)^{s/2}\,dt\right)
\le C v_{n+1}(\eta_A(w)). 
\]
Since $\Omega$ is bounded, the derivatives of $e^{-i\Phi_A}$ through
order $n$ are bounded on $\Omega$. Leibniz' rule, the bounded weak
derivatives of $\sigma$, and the Schwartz decay of $\phi$ and $g$ imply
\[
\sup_{y,w,\beta}\|\widetilde\rho_A(\cdot,y,w,\beta)\|_{W^{n,r}(\Omega)}
\le C|\Omega|^{1/r}. \tag{3.6}\label{eq3.19}
\]
Indeed, derivatives of the activation-window product introduce at most
$v_n(\eta_A(w))$, while its decay in $\beta$ is bounded by a constant
times the second factor of $\vartheta_A$. The fixed shift $t_0$ is
absorbed into this constant.

The weighted modulation-space embedding established above gives
\[
\int_{\mathbb R^2}v_{n+1}(\eta_A(w))|S_g^Af(y,w)|\,dy\,dw
\le C\|f\|_{M_{m_A,A}^{p,q}}.\tag{3.7}\label{eq3.20}
\]
Consequently the normalized-atom representation in \eqref{eq3.14}
is a Bochner integral in $W^{n,r}(\Omega)$ with a finite complex
coefficient measure. (Strong measurability follows from translation
continuity in local $L^r$ for the weak derivatives, with $r<\infty$.)
The representing-measure upper bound therefore implies
\[
\|f\|_{\mathcal K(\widetilde{\mathcal D}_A)}
\le |C_{A,\sigma,\phi}|\int_{\mathbb R^2}|S_g^Af(y,w)|
 \int_{\mathbb R}\vartheta_A(w,\beta)\,d\beta\,dy\,dw
\le C\|f\|_{M_{m_A,A}^{p,q}}.\tag{3.8}\label{eq3.22}
\]
The integral identity initially holds by STFT inversion and hence also
in this Sobolev space by the preceding absolute norm estimate.

The realification of $W^{n,r}(\Omega)$ has type $2$ for $2\le r<\infty$.
Apply Maurey's theorem to the phase-expanded dictionary
$\widetilde{\mathcal D}_A^{\mathrm{ph}}$ and absorb the unit phases into
the complex coefficients. Together with \eqref{eq3.19}, this gives
\[
\inf_{f_N\in\Sigma_N(\widetilde{\mathcal D}_A)}
\|f-f_N\|_{W^{n,r}(\Omega)}
\le C N^{-1/2}|\Omega|^{1/r}\|f\|_{M_{m_A,A}^{p,q}}.
\tag{3.9}\label{eq3.23}
\]
Rescaling the coefficients absorbs $\vartheta_A$, so that
$\Sigma_N(\widetilde{\mathcal D}_A)\subset\Sigma_N(\mathcal D_A)$.
This proves the assertion.
\end{proof}

Define the SAFT by $\mathcal{F}_A f(w) = \int_{\mathbb{R}} f(x) K_A(x,w)\,dx$ and the weighted A-Fourier--Lebesgue space by
\[
\mathcal{F}L_{v_s,A}^{q}(\mathbb{R})
= \bigl\{ f\in \mathscr{S}'(\mathbb{R}): v_s \mathcal{F}_A f \in L^q(\mathbb{R}) \bigr\},\quad
\|f\|_{\mathcal{F}L_{v_s,A}^{q}} = \|v_s \mathcal{F}_A f\|_{L^q}.
\]
The same factorization gives
\[
\|f\|_{\mathcal{F}L_{v_s,A}^{q}} \asymp \|\mathcal{C}_A f\|_{\mathcal{F}L_{v_s}^{q}},\quad
\|h\|_{\mathcal{F}L_{v_s}^{q}} := \bigl\|v_s \widehat{h}\bigr\|_{L^q}.   \tag{3.10}\label{eq3.24}
\]
The fixed scaling of the Fourier variable is absorbed into the equivalence constant.

\begin{corollary}[Local approximation in the weighted A-Feichtinger class]\label{cor:feichtinger}
Assume the hypotheses of Theorem~\ref{thm:main}. If $f\in M^{1,1}_{1\otimes v_{n+1},A}(\R),$
then there is a constant $C>0$, independent of $f$ and $N$, such that
\begin{equation}\label{eq3.25}
 \inf_{f_N\in\Sigma_N(\mathcal{D}_A)} \|f - f_N\|_{W^{n,r}(\Omega)}\le C N^{-1/2}|\Omega|^{1/r}
 \norm{f}_{M^{1,1}_{1\otimes v_{n+1},A}} \tag{3.11}
\end{equation} 
for every $N\in\mathbb N$.
\end{corollary}

\begin{proof}
Choose $p=q=1$ in Theorem~\ref{thm:main}. The admissible spatial weight exponent is $s_1=0$, while the frequency weight exponent is $s_2=n+1$. Thus $m=1\otimes v_{n+1}$, and \eqref{eq3.25} follows directly from \eqref{eq3.10}.
\end{proof}

We use the standard local Fourier--Lebesgue/modulation-space estimate: if $E \subset \mathbb{R}$ is bounded and $u$ is supported in $E$, then, for $0 < q \leq 2$ and $s \in \mathbb{R}$ \cite[Proposition 3.8 and Remark 3.11]{toft2014},
\[
\|u\|_{M_{1\otimes v_s}^{1,q}(\mathbb{R})} \leq C_E |E+E| \|u\|_{\mathcal{F}L_{v_s}^{q}(\mathbb{R})}, \qquad E+E = \{x+y: x,y\in E\}. 
\]
By applying \eqref{eq3.24} to the chirp-modulated function $\mathcal{C}_A u$, and observing that multiplication by $e^{i\Phi_A}$ does not change the support, this transfers to the STSAFT setting:
\[
\|u\|_{M_{1\otimes v_s,A}^{1,q}} \leq C_{A,E} |E+E| \|u\|_{\mathcal{F}L_{v_s,A}^{q}}. \tag{3.12}\label{eq3.27}
\]

\begin{proposition}\label{prop3.3}
(Local STSAFT approximation in weighted A-Fourier--Lebesgue spaces).
Assume the dictionary hypotheses of Theorem~\ref{thm:main}.
Let $0 < q \leq 2$, $2 \leq r < \infty$, and let $\Omega \subset \mathbb{R}$ be bounded.
Suppose
\[
g = \mathbf{1}_\Omega f \in \mathcal{F}L_{v_s,A}^{q}(\mathbb{R}), \qquad
s 
\begin{cases}
=n+1, & 0 < q \leq 1,\\
> n+1 + 1/q', & 1 < q \leq 2.
\end{cases}
\]
Then
\[
\inf_{f_N\in\Sigma_N(\mathcal{D}_A)} \|f - f_N\|_{W^{n,r}(\Omega)}
\leq C N^{-1/2} |\Omega|^{1/r} |\Omega+\Omega| \|\mathbf{1}_\Omega f\|_{\mathcal{F}L_{v_s,A}^{q}}
\label{eq3.28}\tag{3.13}
\]
for all $N\in\mathbb{N}$.
If $\Omega$ is an interval (and hence convex), then $|\Omega+\Omega|=2|\Omega|$, so
\[
\inf_{f_N\in\Sigma_N(\mathcal{D}_A)} \|f - f_N\|_{W^{n,r}(\Omega)}
\leq C N^{-1/2} |\Omega|^{1+1/r} \|\mathbf{1}_\Omega f\|_{\mathcal{F}L_{v_s,A}^{q}}.
\label{eq3.29}\tag{3.14}
\]

\noindent\textit{Proof.}
Set $g=\mathbf{1}_\Omega f$. Since $g=f$ almost everywhere on $\Omega$, the Sobolev approximation error on $\Omega$ is unchanged if $f$ is replaced by $g$. By \eqref{eq3.27} with $E=\Omega$,
\[
\|g\|_{M_{1\otimes v_s,A}^{1,q}} \leq C_{A,\Omega} |\Omega+\Omega| \|g\|_{\mathcal{F}L_{v_s,A}^{q}}.
\]
The selected exponent $s$ is precisely the frequency-weight condition in Theorem~\ref{thm:main} with $p=1$, for which the spatial weight exponent can be chosen as $s_1=0$. Applying Theorem~\ref{thm:main} to $g$ gives
\[
\inf_{f_N\in\Sigma_N(\mathcal{D}_A)} \|g - f_N\|_{W^{n,r}(\Omega)}
\leq C N^{-1/2} |\Omega|^{1/r} \|g\|_{M_{1\otimes v_s,A}^{1,q}}.   
\]
Combining these two inequalities yields \eqref{eq3.28}. If $\Omega$ is an interval, then $\Omega+\Omega$ is an interval of twice the length, and hence $|\Omega+\Omega|=2|\Omega|$, proving \eqref{eq3.29}.
\qed
\end{proposition}

\begin{definition}
(A-Barron space). For $s\in\mathbb{R}$, define
\[
\mathcal{B}_A^s(\mathbb{R}) := \mathcal{F}L_{v_s,A}^{1}(\mathbb{R}), \qquad
\|f\|_{\mathcal{B}_A^s} := \|f\|_{\mathcal{F}L_{v_s,A}^{1}}. 
\]
\end{definition}

\begin{corollary}\label{coro3.5}
(Local STSAFT approximation in A-Barron spaces).
Assume the hypotheses of Theorem~\ref{thm:main}, and let $\Omega\subset\mathbb{R}$ be a bounded interval. If
$g = \mathbf{1}_\Omega f \in \mathcal{B}_A^{n+1}(\mathbb{R})$
then
\[
\inf_{f_N\in\Sigma_N(\mathcal{D}_A)} \|f-f_N\|_{W^{n,r}(\Omega)}
\leq C N^{-1/2} |\Omega|^{1/r} |\Omega+\Omega| \|\mathbf{1}_\Omega f\|_{\mathcal{B}_A^{n+1}}.
\label{eq3.32}\tag{3.15}
\]
In particular, since $|\Omega+\Omega|=2|\Omega|$ for an interval,
\[
\inf_{f_N\in\Sigma_N(\mathcal{D}_A)} \|f-f_N\|_{W^{n,r}(\Omega)}
\leq C N^{-1/2} |\Omega|^{1+1/r} \|\mathbf{1}_\Omega f\|_{\mathcal{B}_A^{n+1}}.
\label{eq3.33}
\tag{3.16}
\]
If $\operatorname{supp}f \subseteq \Omega$, the cutoff $\mathbf{1}_\Omega$ may be omitted.

\noindent\textit{Proof.} By definition, $\mathcal{B}_A^{n+1} = \mathcal{F}L_{v_{n+1},A}^{1}$. Thus the present corollary is Proposition~\ref{prop3.3} with $q=1$ and $s=n+1$. The interval estimate follows from $|\Omega+\Omega|=2|\Omega|$.
\qed
\end{corollary}

\begin{remark}[Scope of the hard-cutoff hypothesis]
The assumptions in Proposition~\ref{prop3.3} and Corollary~\ref{coro3.5}
are imposed on the zero extension $\mathbf 1_\Omega f$, not merely on the
restriction of $f$ to $\Omega$. In particular, membership of $f$ in a
global Fourier--Lebesgue or A-Barron space does not by itself imply that
$\mathbf 1_\Omega f$ belongs to the same space: a hard cutoff may create
jumps at the boundary. Thus these statements are conditional local
estimates and are not asserted for every globally regular $f$.
A useful alternative is to choose $\chi\in C_c^\infty(\mathbb R)$ with
$\chi=1$ on a neighbourhood of $\overline\Omega$ and assume directly that
$\chi f$ belongs to the Fourier--Lebesgue class required in
Proposition~\ref{prop3.3}. We do not infer this condition from multiplication
by $\chi$ for $q<1$ without a separate quasi-Banach multiplier argument.
The constants in the local embedding and approximation estimate may depend
on $A$, the fixed window, and the support set (hence on $\Omega$); the
displayed powers of $|\Omega|$ are not claimed to be domain-uniform gains.
\end{remark}

\subsection{Global Approximation}


\begin{proposition}[Restricted-translation obstruction]
  For every bounded measurable \(\Omega\subset\R\) and \(h\in\Sclass(\R)\),
  \begin{equation}
    \int_\Omega\int_\R V_gh(y,\xi)g(x-y)e^{2\pi i\xi x}\,d\xi\,dy
    =m_\Omega(x)h(x),\qquad
    m_\Omega(x)=\int_\Omega |g(x-y)|^2\,dy.
    \label{eq:truncation}\tag{3.17}
  \end{equation}
  In particular, the truncated integral is not an inversion formula in general.
  Moreover, a dictionary all of whose atoms have the factor \(g(x-y)\), \(y\in\Omega\),
  cannot approximate every Schwartz function globally when \(g\) is compactly supported.
\end{proposition}

\begin{proof}
  By Fourier inversion applied to \(h\,\overline{g(\cdot-y)}\), we have
  \(\int_\R V_gh(y,\xi)e^{2\pi i\xi x}\,d\xi=h(x)\overline{g(x-y)}\).
  By multiplying by \(g(x-y)\) and integrating in \(y\), we obtain \eqref{eq:truncation}.
  For a Gaussian window and \(|\Omega|>0\), we have \(0<m_\Omega(x)<1\) for every \(x\).
  Thus even for this strictly positive window we do not recover \(h\).

  Now choose \(g\in C_c^\infty(\R)\) with \(\|g\|_2=1\), and put
  \(K=\overline\Omega+\supp g\).
  By taking the activation \(\sigma\equiv1\) and the activation window
  \(\varphi(t)=e^{-\pi t^2}\), with \(t_0=0\) and any \(\tau_0\ne0\),
  We satisfy the previous remaining activation and nonvanishing assumptions.
  Since every such atom is supported in the compact set \(K\), we have
  \(\supp v\subset K\) for every finite network \(v\).
  For a nonzero \(h\in C_c^\infty(\R\setminus K)\), we therefore have
  \[\inf_v\|h-v\|_{W^{n,r}(\R)}\ge\|h\|_{L^r(\R)}>0.\]
  Since \(h\) belongs to every polynomially weighted modulation space occurring in the stated
  global estimate, we cannot have an upper bound tending to zero as \(N\to\infty\).
\end{proof}

The Gaussian computation identifies a failure of the displayed reconstruction, rather than
a non-density result for every Gaussian restricted-shift dictionary.
The repair below changes the dictionary to all translations and imposes an explicit spatial moment.

Write
\(A=\left(\begin{smallmatrix}a_A&b_A&p_A\\c_A&d_A&q_A\end{smallmatrix}\right)\),
where \(b_A\ne0\) and \(a_Ad_A-b_Ac_A=1\), and put
\[\Phi_A(x)=\frac{a_A}{2b_A}x^2+\frac{p_A}{b_A}x, \quad U_A f=e^{i\Phi_A}f, \quad \Psi_A(w)=\frac{d_Aw^2}{2b_A}+\left(q_A-\frac{d_Ap_A}{b_A}\right)w.\]
We have
\[S_g^Af(y,w)=\frac{e^{i\Psi_A(w)}}{\sqrt{2\pi|b_A|}} V_g(U_Af)\left(y,\frac{w}{2\pi b_A}\right).\]
For an integer \(n\ge1\) and \(2\le r<\infty\), define
\[\J_{A,n,r}(f)=\iint_{\R^2}|V_g(U_Af)(y,\xi)| \br{y}^{n}\br{\xi}^{n+1-1/r}\,dy\,d\xi. \tag{3.18}\]
Equivalently, by the change of variables \(w=2\pi b_A\xi\),
\[\J_{A,n,r}(f)=\frac{1}{\sqrt{2\pi|b_A|}} \iint_{\R^2}|S_g^Af(y,w)|\br{y}^{n} \left\langle\frac{w}{2\pi b_A}\right\rangle^{n+1-1/r}\,dy\,dw.\]
Thus this is a weighted STSAFT coefficient norm with the frequency rescaling made explicit.

\begin{theorem}[Global approximation with all translations]\label{thm:global}
  Let \(g\in\Sclass(\R)\) be real with \(\|g\|_2=1\).
  Let \(H\colon\R\to\R\) be nonzero, belong to \(W^{n,\infty}(\R)\), and satisfy
  \[\|\br{t}^{M}H^{(j)}(t)\|_\infty<\infty \qquad(M\ge0,\ 0\le j\le n),\]
  where derivatives are weak derivatives. Choose \(\tau\ne0\) with \(\widehat H(\tau)\ne0\), and set
  \[\rho_A(x;y,\xi,b)=e^{-i\Phi_A(x)}H(\xi x/\tau+b)g(x-y), \qquad \D_A=\{\rho_A(\cdot;y,\xi,b):y,\xi,b\in\R\}.\]
  Let \(\Sigma_N(\D_A)\) denote the complex linear combinations of at most \(N\) atoms.
  If \(f\in L^2(\R)\) and \(\J_{A,n,r}(f)<\infty\), then \(f\in W^{n,r}(\R)\) and
  \begin{equation}
    \inf_{v\in\Sigma_N(\D_A)}\|f-v\|_{W^{n,r}(\R)}
    \le C N^{-1/2}\J_{A,n,r}(f),\qquad N\ge1.
    \label{eq:global}\tag{3.19}
  \end{equation}
  Here \(C\) depends only on \(A,n,r,g,H,\tau\).
  If \(a_A=0\), the same statement holds with the spatial weight \(\br{y}^{n}\) omitted.
\end{theorem}

\begin{proof}
  We use \(\|v\|_{W^{n,r}}=(\sum_{j=0}^n\|v^{(j)}\|_r^r)^{1/r}\).
  Since \(H\in L^1\) is nonzero, by continuity and injectivity of the Fourier transform
  we can choose the stated nonzero \(\tau\).

  Put \(c=\xi/\tau\), \(L=\br{\xi}\), and \(r'=r/(r-1)\).
  For measurable \(h,\psi\) with \(\br{\cdot}h,\br{\cdot}\psi\in L^r\), put
  \(F(\beta)=\|h(cu+\beta)\psi(u)\|_{L^r_u}\).
  By H\"older's inequality in \(\beta\) and Tonelli's theorem, we have
  \[
    \begin{aligned}
      \int_\R F(\beta)\,d\beta
      &\le
      \left(\int_\R\br{\beta/L}^{-r'}\,d\beta\right)^{1/r'}
      \left(\iint_{\R^2}\br{\beta/L}^{r}|h(cu+\beta)\psi(u)|^r\,du\,d\beta\right)^{1/r}.
    \end{aligned}
  \]
  Since \(r'>1\), the first factor equals \(C_rL^{1-1/r}\).
  By setting \(z=cu+\beta\) in the second factor and using
  \(\br{(z-cu)/L}\le C_\tau\br{z}\br{u}\), we obtain
  \[\int_\R\|h(cu+\beta)\psi(u)\|_{L^r_u}\,d\beta \le C_{r,\tau}\br{\xi}^{1-1/r} \|\br{z}h(z)\|_r\|\br{u}\psi(u)\|_r.\]
  Now put \(u=x-y\), \(\beta=b+\xi y/\tau\), and
  \(R_A(x;y,\xi,\beta)=e^{-i\Phi_A(x)}H(\xi(x-y)/\tau+\beta)g(x-y)\).
  Since \(\Phi_A\) is quadratic, we have
  \[\partial_x^j e^{-i\Phi_A(x)}=P_{j,A}(x)e^{-i\Phi_A(x)},\qquad |P_{j,A}(y+u)|\le C_{A,j}\br{y}^{j}\br{u}^{j}.\]
  By the weak Leibniz and chain rules, every derivative of \(R_A\) of order at most \(n\)
  is a finite sum of terms with absolute value bounded by
  \(C\br{y}^j\br{u}^j|c|^\ell|H^{(\ell)}(cu+\beta)g^{(k)}(u)|\),
  where \(j+\ell+k\le n\).
  By the preceding estimate with \(h=H^{(\ell)}\) and
  \(\psi(u)=\br{u}^jg^{(k)}(u)\), we obtain
  \begin{equation}
    \int_\R\|R_A(\cdot;y,\xi,\beta)\|_{W^{n,r}}\,d\beta
    \le C\br{y}^{n}\br{\xi}^{n+1-1/r}.
    \label{eq:atoms}\tag{3.20}
  \end{equation}
  Since \(H\) has a continuous representative, we can also define the atoms at \(\xi=0\);
  the positive-order derivatives of \(H(\xi(x-y)/\tau+\beta)\) then vanish.
  By choosing Borel representatives of the weak derivatives, we have jointly measurable
  derivative expressions. By separability of the finite product of \(L^r\) spaces,
  we obtain strong measurability of the atom map into \(W^{n,r}\).

  By substituting \(t=\xi x/\tau+b\) in the definition of \(\widehat H(\tau)\), we have
  \[e^{2\pi i\xi x}=\frac{1}{\widehat H(\tau)} \int_\R H(\xi x/\tau+b)e^{-2\pi i\tau b}\,db.\]
  Define
  \[Q_f(y,\xi,\beta)=\frac{V_g(U_Af)(y,\xi)}{\widehat H(\tau)} e^{2\pi i\xi y}e^{-2\pi i\tau\beta}.\]
  By \eqref{eq:atoms}, we have
  \[M_f:=\iiint_{\R^3}|Q_f(y,\xi,\beta)| \|R_A(\cdot;y,\xi,\beta)\|_{W^{n,r}}\,d\beta\,dy\,d\xi \le C\J_{A,n,r}(f).\]
  Thus the integral \(F=\iiint Q_f R_A\) converges in \(W^{n,r}\).
  By the integrated atom estimate with no derivatives and exponent \(2\), we also have
  \(\int_\R\|R_A\|_2\,d\beta\le C\br{\xi}^{1/2}\).
  Since \(n+1-1/r\ge1/2\), the same integral converges in \(L^2\),
  and both limits agree as distributions.
  By Bochner Fubini, the substitution \(b=\beta-\xi y/\tau\), and the scalar identity above,
  we obtain in \(L^2\)
  \[F=U_A^{-1}\iint_{\R^2}V_g(U_Af)(y,\xi)g(\cdot-y)e^{2\pi i\xi\,\cdot}\,dy\,d\xi=f.\]
  By the standard STFT inversion formula \cite{Grochenig2001}, we have the last equality;
  since \(V_g(U_Af)\in L^1\), the integral is absolutely convergent in \(L^2\).
  Therefore \(f\in W^{n,r}\), with the required integral representation.

  Put \(X=W^{n,r}(\R;\C)\).
  Since \(2\le r<\infty\), we can regard \(X\) as a real Banach space of type 2,
  with constant \(T_2(X)\).
  By the standard Maurey sampling argument \cite{Pisier,AbdeljawadCarioniCordero2026}, we can discretise the preceding
  representation as follows.
  If \(M_f=0\), then \(f=0\), so assume \(M_f>0\).
  On the set where \(|Q_f|\|R_A\|_X>0\), define
  \[d\mathbb P=\frac{|Q_f|\|R_A\|_X}{M_f}\,d\beta\,dy\,d\xi, \qquad Z=\frac{Q_f}{|Q_f|}\frac{R_A}{\|R_A\|_X}.\]
  Then \(\|Z\|_X=1\) and \(\mathbb EZ=f/M_f\).
  For independent copies \(Z_1,\ldots,Z_N\), their independent copies \(Z'_1,\ldots,Z'_N\),
  and independent symmetric signs \(\varepsilon_j\), by symmetrisation and type 2 we have
  \[
    \begin{aligned}
      \mathbb E\left\|\sum_{j=1}^N(Z_j-\mathbb EZ_j)\right\|_X
      &\le\mathbb E\left\|\sum_{j=1}^N(Z_j-Z'_j)\right\|_X
      \le2\mathbb E\left\|\sum_{j=1}^N\varepsilon_jZ_j\right\|_X
      \le2T_2(X)\sqrt N.
    \end{aligned}
  \]
  Since every realisation of \((M_f/N)\sum_jZ_j\) belongs to \(\Sigma_N(\D_A)\), we obtain
  \(\inf_{v\in\Sigma_N(\D_A)}\|f-v\|_X\le2T_2(X)M_fN^{-1/2}\).
  Thus \eqref{eq:global} follows.
  If \(a_A=0\), all \(P_{j,A}\) are constant, so we can omit \(\br{y}^{n}\)
  in \eqref{eq:atoms} and throughout the proof.
\end{proof}

The proof uses standard STFT inversion and type-2 sampling, including the formulation in
\cite[Proposition 6]{AbdeljawadCarioniCordero2026}. The dictionary estimate \eqref{eq:atoms} supplies the required
global integrability.

\begin{corollary}[Gaussian-window ReLU and real-valued networks]\label{cor:global-relu}
  Fix \(t_0\in\R\), \(\tau_0\ne0\), and the SAFT parameters \(A\).
  Let \(g(x)=\pi^{-1/4}e^{-x^2/2}\), and put
  \[h_j(x)=\left(\frac{\eta_jx}{\tau_0}+\beta_j\right)_+ \exp\left[-\frac{1}{2}\left(\frac{\eta_jx}{\tau_0}+\beta_j-t_0\right)^2\right] e^{-(x-y_j)^2/2}.\]
  Let \(\mathcal R_{A,N}\) consist of the real-valued networks
  \[v_N(x)=\sum_{j=1}^N \big[\alpha_jh_j(x)\cos\Phi_A(x)+\gamma_jh_j(x)\sin\Phi_A(x)\big], \qquad y_j,\eta_j,\beta_j,\alpha_j,\gamma_j\in\R.\]
  For every real \(f\in L^2(\R)\) with \(\J_{A,1,r}(f)<\infty\),
  \begin{equation}
    \inf_{v_N\in\mathcal R_{A,N}}\|f-v_N\|_{W^{1,r}(\R)}
    \le C N^{-1/2}\iint_{\R^2}|V_g(U_Af)(y,\xi)|
    \br{y}\br{\xi}^{2-1/r}\,dy\,d\xi.
    \label{eq:relu} \tag{3.21}
  \end{equation}
  In particular, the target \(f_{\omega,\kappa}(x)=e^{-x^2}\sin(\omega x+\pi\kappa x^2)\)
  in the experiments satisfies this estimate for every fixed \(\omega,\kappa\in\R\).
  The same upper bound holds after restriction to \([-1,1]\), and after allowing an additional
  constant output bias on that interval.
\end{corollary}

\begin{proof}
  Put \(H(t)=t_+e^{-(t-t_0)^2/2}\).
  Since \(H\) is continuous at zero, its weak first derivative is the function
  \[H'(t)=\mathbf1_{(0,\infty)}(t)e^{-(t-t_0)^2/2} -t_+(t-t_0)e^{-(t-t_0)^2/2}.\]
  Thus \(H\in W^{1,\infty}\), and every polynomially weighted \(H\) and \(H'\) is bounded.
  Since \(H\ge0\) is nonzero, we have \(\widehat H(0)=\|H\|_1>0\);
  by continuity we can choose \(\tau\ne0\) with \(\widehat H(\tau)\ne0\).
  By applying the theorem with \(n=1\), we obtain the complex-network estimate.
  Since \(\eta_j\) is unrestricted, the substitution \(\eta_j=\tau_0\xi_j/\tau\)
  preserves the dictionary for any prescribed \(\tau_0\ne0\).
  The normalising factor \(\pi^{-1/4}\) is absorbed into the output coefficients.
  For \(c_j=\alpha_j+i\gamma_j\), we have
  \(\Repart(c_jh_je^{-i\Phi_A})=\alpha_jh_j\cos\Phi_A+\gamma_jh_j\sin\Phi_A\).
  Since taking real parts is contractive in \(W^{1,r}\), we obtain \eqref{eq:relu}.
  Since \(U_Af_{\omega,\kappa}\in\Sclass(\R)\), its STFT is Schwartz and the coefficient
  integral is finite. By restricting the error norm to \([-1,1]\), we obtain the local estimate.
  Since the class with an output bias contains the class with zero bias, we obtain the last assertion.
\end{proof}

The Gaussian window in Corollary~\ref{cor:global-relu} is fixed there as
$g(x)=\pi^{-1/4}e^{-x^2/2}$. It is distinct from the phase-analysis
window $g_{\mathrm{ph}}(x)=2^{1/4}e^{-\pi x^2}$ used in
Section~\ref{sec:pm}; the two normalizations are not being identified.



\paragraph{A computable power proxy.}
For the strengthened phase estimate, the dimensionless proxy is
\[
Q_\kappa^{\mathrm{pow}}(\lambda)
=\frac{\langle\lambda\rangle^n
\langle\kappa-\lambda\rangle^{n+d_*}}
{\langle\kappa\rangle^{n+d_*}},
\qquad d_*=\frac32-\frac1r.
\]
If $\lambda_\kappa=\kappa-\kappa^\alpha$ with $0\le\alpha<1$, then
$Q_\kappa^{\mathrm{pow}}(\lambda_\kappa)
\asymp\kappa^{(n+d_*)\alpha-d_*}$.
It tends to zero precisely when
$\alpha<d_*/(n+d_*)$. This is a proxy for the ratio of the specified
upper-bound factors. It is not the exact ratio of best approximation
errors. In particular, for the experimental choice $n=1$, $r=2$,
we have $d_*=1$ and the exponent is $2\alpha-1$.

\subsection{Phase matching and the size of the approximation bound}
\label{sec:pm}
Theorem~\ref{thm:STSAFT-local-approximation} supplies the approximation rate.
We now specialize its dictionary to a fixed-scale SAFT family to quantify
how phase matching changes the representation bound. The additional Gaussian
window and envelope assumptions in this subsection do not replace the more
general hypotheses of that theorem.

Let $\Omega\subset\mathbb R$ be a bounded interval of positive length,
$n\ge1$ and $2\le r<\infty$. Put $s_*:=n+1-\frac1r$ and $d_*:=\frac32-\frac1r$. Then $s_*+\frac12=n+d_*$. Write
$\langle t\rangle=(1+t^2)^{1/2}$ and fix
$g_{\mathrm{ph}}(x)=2^{1/4}e^{-\pi x^2}$, so that $\|g_{\mathrm{ph}}\|_2=1$.
For $t\ge0$, define $\mathcal M_t$ and $\mathcal A_t$ as the spaces of
complex-valued $L^2(\mathbb R)$ functions with finite respective norms
\begin{align}
 \|h\|_{\mathcal M_t}
 &=\iint_{\mathbb R^2}|V_{g_{\mathrm{ph}}}h(y,\xi)|\langle\xi\rangle^t\,dy\,d\xi,
 \label{eq:pm-Mnorm}\tag{3.22}\\
 \|h\|_{\mathcal A_t}
 &=\iint_{\mathbb R^2}|V_{g_{\mathrm{ph}}}h(y,\xi)|\langle y\rangle^t
      \langle\xi\rangle^t\,dy\,d\xi.\label{eq:pm-Anorm}\tag{3.23}
\end{align}
Thus these are the fixed-window weighted $M^{1,1}$ norms already used above,
with frequency weight and product weight, respectively. In the phase
analysis below we use the non-integer exponent $s_*$; the definitions make
sense for every real $t\ge0$.
Let $B=W^{n,r}(\Omega;\mathbb C)$, with weak-derivative norm
\[
 \|h\|_B=\left(\sum_{j=0}^n\|h^{(j)}\|_{L^r(\Omega)}^r\right)^{1/r},
 \qquad
 E_N(h;\mathcal D)=\inf_{v\in\Sigma_N(\mathcal D)}\|h-v\|_B.
\]
Retain a fixed real-valued activation $\sigma\in W^{k,\infty}(\mathbb R)\setminus\{0\}$,
$k\ge n$, and a fixed real-valued Schwartz window $\phi$, with
$t_0,\tau_0$ satisfying the assumptions of
Theorem~\ref{thm:STSAFT-local-approximation}. Set
\begin{equation}\label{eq:pm-dictionary}
 \mathcal D_0=\left\{
 \sigma(\xi x/\tau_0+\beta)
 \phi(\xi x/\tau_0+\beta-t_0)g_{\mathrm{ph}}(x-y):y,\xi,\beta\in\mathbb R
 \right\}\bigg|_\Omega.\tag{3.24}
\end{equation}
Define $C_\lambda h(x)=e^{i\pi\lambda x^2}h(x)$ and consider
\begin{equation}\label{eq:pm-SAFT}
 A_\lambda=\begin{bmatrix}-\lambda&1/(2\pi)&0\\-2\pi&0&0\end{bmatrix},
 \qquad \mathcal D_\lambda=C_\lambda\mathcal D_0.\tag{3.25}
\end{equation}
This matrix has determinant one in its linear part. The kernel factorization
gives $\Phi_{A_\lambda}(x)=-\pi\lambda x^2$, $\Psi_{A_\lambda}=0$,
and $\eta_{A_\lambda}(w)=w$, hence
 $S_g^{A_\lambda}f=V_{g_{\mathrm{ph}}}(C_{-\lambda}f),
 \mathcal D_{A_\lambda}=\mathcal D_\lambda.$
Theorem~\ref{thm:global} below gives a stronger input for this same
dictionary. Indeed, set
\[
H(t)=\sigma(t)\phi(t-t_0).
\]
The assumptions on $\sigma$ and $\phi$ imply that $H$ and its weak
derivatives through order $n$ have arbitrary polynomial decay. Moreover,
$\widehat H(\tau_0)=V_\phi\sigma(t_0,\tau_0)\ne0$ (with the Fourier
convention fixed above). Apply Theorem~\ref{thm:global} with $A=A_0$ and window $g=g_{\mathrm{ph}}$, and
then restrict the resulting global Sobolev error to $\Omega$. After the
harmless reparametrization of the unrestricted frequency variable, its
dictionary is precisely $\mathcal D_0$. Consequently,
\begin{equation}
 E_N(h;\mathcal D_0)\le K_*N^{-1/2}\|h\|_{\mathcal M_{s_*}},
 \qquad h\in\mathcal M_{s_*},\quad N\ge1.
 \tag{H$^*$}\label{eq:Hstar}
\end{equation}
Here $K_*$ depends on the fixed dictionary data and $\Omega$, but not on
$h,N,\kappa,\lambda$. This is stronger than the estimate obtained from
the local theorem with frequency weight $n+1$.
The activation, both windows
and the frequency scale remain fixed as $\lambda$ varies.
For $f_\kappa=C_\kappa u$, the transformed envelope is
$C_{\kappa-\lambda}u$, so $\lambda=\kappa$ cancels its quadratic phase.
All constants below may depend on the fixed data, but not on the varying
phase parameters or $N$, unless explicitly indexed by an envelope.

\begin{theorem}
  Assume \eqref{eq:Hstar}, and let \(f_\kappa=C_\kappa a\) and \(\D_\lambda=C_\lambda\D_0\).
  There is a constant \(C\), independent of \(a,\kappa,\lambda,N\), such that
  \begin{equation}
    E_N(f_\kappa;\D_\lambda)
    \le C N^{-1/2}\br{\lambda}^{n}\br{\kappa-\lambda}^{n+d_*}
    \|a\|_{\mathcal A_{s_*}}
    \qquad(a\in\mathcal A_{s_*}).
    \label{eq:main}\tag{3.26}
  \end{equation}
  The constant may depend on \(K,\Omega,n,r\).

  Moreover, for every fixed \(0\ne a\in\Sclass(\R)\),
  \begin{equation}
    \lim_{|\delta|\to\infty}
    \frac{\|C_\delta a\|_{\mathcal M_{s_*}}}{|\delta|^{s_*+1/2}}
    =\|g_{\mathrm{ph}}\|_1\int_\R |t|^{s_*}|a(t)|\,dt>0.
    \label{eq:asym}\tag{3.27}
  \end{equation}

  For such an \(a\), define
  \[\widetilde B_\lambda(\kappa,N)=K_* N^{-1/2} \|C_\lambda\|_{W^{n,r}(\Omega)\to W^{n,r}(\Omega)} \|C_{\kappa-\lambda}a\|_{\mathcal M_{s_*}}.\]
  Then \(E_N(f_\kappa;\D_\lambda)\le \widetilde B_\lambda(\kappa,N)\), and, as \(|\kappa|\to\infty\),
  \[\widetilde B_0\asymp_a N^{-1/2}|\kappa|^{n+d_*},\qquad \widetilde B_\kappa\lesssim_a N^{-1/2}\br{\kappa}^{n},\qquad \frac{\widetilde B_\kappa}{\widetilde B_0}=O_a(|\kappa|^{-d_*}).  \tag{3.28}\]
\end{theorem}

\begin{proof}
  By evaluating the Gaussian integral, we have
  \[V_{g_{\mathrm{ph}}}(C_\delta g_{\mathrm{ph}})(y,\xi)=\frac{\sqrt2}{\sqrt{2-i\delta}} \exp\!\left(-\pi y^2+\frac{\pi(y-i\xi)^2}{2-i\delta}\right).\]
  By taking the modulus and integrating in \(y\), we obtain
  \[q_\delta(\xi):=\int_\R|V_{g_{\mathrm{ph}}}(C_\delta g_{\mathrm{ph}})(y,\xi)|\,dy =\frac{A_\delta}{L_\delta}e^{-\pi\xi^2/L_\delta^2},\qquad A_\delta=\sqrt2(4+\delta^2)^{1/4},\quad L_\delta=\sqrt{2+\delta^2}.\]
  By substituting \(\xi=L_\delta u\) and using the Gaussian moments, we have, for every \(v\in\R\),
  \begin{align*}
    \int_\R\br{\xi+v}^{s_*} q_\delta(\xi)\,d\xi
    &=A_\delta\int_\R\br{L_\delta u+v}^s e^{-\pi u^2}\,du\\
    &\le C_s\br{\delta}^{1/2}\bigl(\br{v}^{s_*}+\br{\delta}^{s_*}\bigr).
  \end{align*}
  Since \(C_\delta M_\eta T_yg_{\mathrm{ph}}=e^{-i\pi\delta y^2}M_{\eta+\delta y}T_y(C_\delta g_{\mathrm{ph}})\),
  by STFT covariance we identify its \(\mathcal M_{s_*}\)-norm with this integral for \(v=\eta+\delta y\).
  Since \(\br{\eta+\delta y}\le C\br{\delta}\br{y}\br{\eta}\), we get
  \(\|C_\delta M_\eta T_yg_{\mathrm{ph}}\|_{\mathcal M_{s_*}}\le C_s\br{\delta}^{s_*+1/2}\br{y}^{s_*}\br{\eta}^{s_*}\).

  Since \(a\in\mathcal A_{s_*}\) and \(\|g_{\mathrm{ph}}\|_2=1\), by Gabor reconstruction we have
  \(a=\iint V_{g_{\mathrm{ph}}}a(y,\eta)M_\eta T_yg_{\mathrm{ph}}\,dy\,d\eta\).
  Since the integrand has \(L^2\)-norm \(|V_{g_{\mathrm{ph}}}a(y,\eta)|\in L^1(\R^2)\), we have convergence in \(L^2\).
  By pairing with an arbitrary \(L^2\) function and using STFT orthogonality, we obtain the stated equality.
  Since \(C_\delta\) is unitary and each STFT coefficient is a continuous linear functional on \(L^2\),
  we may apply \(C_\delta\) and pass the STFT through the integral.
  By the triangle inequality and Tonelli's theorem, we obtain
  \[
  \|C_\delta a\|_{\mathcal M_{s_*}}
  \le \iint |V_{g_{\mathrm{ph}}}a(y,\eta)|
  \|C_\delta M_\eta T_yg_{\mathrm{ph}}\|_{\mathcal M_{s_*}}\,dy\,d\eta
  \le C_s\br{\delta}^{s_*+1/2}\|a\|_{\mathcal A_{s_*}}.
  \]
  Moreover, since the derivatives through order \(n\) of \(M_\eta T_yg_{\mathrm{ph}}\) are bounded by
  \(C_n\br{\eta}^n\) uniformly in \(x,y\), by weighted reconstruction we have
  \(\mathcal M_{s_*}\subseteq W^{n,r}(\Omega)\).

  By writing \(\partial_x^j e^{i\pi\lambda x^2}=P_{j,\lambda}(x)e^{i\pi\lambda x^2}\),
  we have \(P_{0,\lambda}=1\) and
  \(P_{j+1,\lambda}=P'_{j,\lambda}+2\pi i\lambda xP_{j,\lambda}\).
  Since \(\Omega\) is bounded, we obtain
  \(\sup_{x\in\Omega}|P_{j,\lambda}(x)|\le C_{\Omega,j}\br{\lambda}^j\) for \(0\le j\le n\).
  By Leibniz's formula for weak derivatives, we have
  \(\|C_\lambda h\|_{W^{n,r}(\Omega)}\le C_{\Omega,n,r}\br{\lambda}^n\|h\|_{W^{n,r}(\Omega)}\).
  Since \(f_\kappa=C_\lambda(C_{\kappa-\lambda}a)\) and
  \(\Sigma_N(\D_\lambda)=C_\lambda\Sigma_N(\D_0)\), we obtain
  \[E_N(f_\kappa;\D_\lambda) \le\|C_\lambda\|_{W^{n,r}\to W^{n,r}} E_N(C_{\kappa-\lambda}a;\D_0) \le \widetilde B_\lambda(\kappa,N).\]
  Therefore, by \eqref{eq:Hstar} and the preceding norm estimates, we prove \eqref{eq:main}.

  Moreover, for the asymptotic formula, we assume \(0\ne a\in\Sclass(\R)\) and set
  \(b(y,x)=a(x)\overline{g_{\mathrm{ph}}(x-y)}\in\Sclass(\R^2)\).
  By completing the quadratic phase, we have
  \[V_{g_{\mathrm{ph}}}(C_\delta a)(y,\delta t) =e^{-i\pi\delta t^2}\int_\R b(y,x)e^{i\pi\delta(x-t)^2}\,dx.\]
  By the Fourier convention \(\widehat u(\eta)=\int u(x)e^{-2\pi ix\eta}\,dx\),
  we write the Fresnel identity in \(\Sclass'(\R)\) as
  \[\mathcal F\!\left(e^{i\pi\delta x^2}\right)(\eta) =|\delta|^{-1/2}e^{i\operatorname{sgn}(\delta)\pi/4}e^{-i\pi\eta^2/\delta},\qquad \delta\ne0.\]
  By applying the Gaussian Fourier transform formula with \(\varepsilon>0\), we have
  \[\mathcal F\!\left(e^{-\pi(\varepsilon-i\delta)x^2}\right)(\eta) =(\varepsilon-i\delta)^{-1/2} \exp\!\left(-\frac{\pi\eta^2}{\varepsilon-i\delta}\right),\]
  where we use the principal square root.
  By analytic continuation of the real Gaussian formula to \(\operatorname{Re}z>0\),
  we obtain this complex Gaussian formula.
  Since \(|e^{-\pi(\varepsilon-i\delta)x^2}|\le1\), by dominated convergence against Schwartz test functions
  we have \(e^{-\pi(\varepsilon-i\delta)x^2}\to e^{i\pi\delta x^2}\) in \(\Sclass'\) as \(\varepsilon\downarrow0\).
  Moreover, since
  \((\varepsilon-i\delta)^{-1/2}\to|\delta|^{-1/2}e^{i\operatorname{sgn}(\delta)\pi/4}\)
  and \(-\pi\eta^2/(\varepsilon-i\delta)\to-i\pi\eta^2/\delta\),
  we have pointwise convergence of the transformed Gaussian to the asserted Fresnel expression.
  Since
  \[\left|(\varepsilon-i\delta)^{-1/2} e^{-\pi\eta^2/(\varepsilon-i\delta)}\right| =(\varepsilon^2+\delta^2)^{-1/4} e^{-\pi\varepsilon\eta^2/(\varepsilon^2+\delta^2)} \le |\delta|^{-1/2},\]
  by dominated convergence we also have convergence in \(\Sclass'\).
  Thus, by continuity of the Fourier transform on \(\Sclass'\), we obtain the Fresnel identity.
  By convolving with \(b(y,\cdot)\in\Sclass(\R)\), we obtain
  \[|\delta|^{1/2}e^{i\pi\delta t^2}V_{g_{\mathrm{ph}}}(C_\delta a)(y,\delta t) =e^{i\operatorname{sgn}(\delta)\pi/4}(U_\delta b)(y,t),\]
  where we set
  \(U_\delta b=\mathcal F^{-1}_{\eta\to t} (e^{-i\pi\eta^2/\delta}\mathcal F_{x\to\eta}b)\).

  By setting \(m_\delta(\eta)=e^{-i\pi\eta^2/\delta}\), differentiating and using \(|e^{iu}-1|\le|u|\),
  we obtain
  \(|\partial_\eta^j(m_\delta-1)(\eta)|\le C_j|\delta|^{-1}\br{\eta}^{2j+2}\)
  for \(|\delta|\ge1\) and every integer \(j\ge0\).
  Since \(\mathcal F_xb\in\Sclass(\R^2)\), by Leibniz's formula we have
  \((m_\delta-1)\mathcal F_xb\to0\) in \(\Sclass(\R^2)\).
  Thus, by continuity of the partial Fourier transform, we have \(U_\delta b\to b\) in \(\Sclass(\R^2)\).
  Moreover, for every \(M>0\), we have
  \(|(U_\delta b)(y,t)|\le C_M(1+y^2+t^2)^{-M}\), uniformly for \(|\delta|\ge1\).
  By substituting \(\xi=\delta t\), we obtain
  \[\frac{\|C_\delta a\|_{\mathcal M_{s_*}}}{|\delta|^{s_*+1/2}} =\iint(\delta^{-2}+t^2)^{s_*/2}|(U_\delta b)(y,t)|\,dy\,dt.\]
  Since the integrand converges pointwise to \(|t|^s|b(y,t)|\) and is bounded by
  \(C_M\br{t}^{s_*}(1+y^2+t^2)^{-M}\), which is integrable for \(2M>s_*+2\),
  by dominated convergence and integration in \(y\) we obtain \eqref{eq:asym}.
  Since \(a\) is nonzero and continuous, we have \(\int|t|^{s_*}|a(t)|\,dt>0\).
  Therefore, by \(C_0=I\), \eqref{eq:asym} and the Sobolev multiplier estimate,
  we obtain the stated comparison of upper bounds.
\end{proof}

\begin{corollary}
  Fix \(0\ne a\in\Sclass(\R)\) and retain all other parameters. For all \(\kappa,\lambda\in\R\) and \(N\ge1\),
  \[\frac{\widetilde B_\lambda(\kappa,N)}{\widetilde B_0(\kappa,N)} \asymp_a \frac{\br{\lambda}^{n}\br{\kappa-\lambda}^{n+d_*}}{\br{\kappa}^{n+d_*}}, \tag{3.29}\]
  with comparison constants independent of \(\kappa,\lambda,N\).
  For any choice \(\lambda_\kappa\in\R\), as \(\kappa\to+\infty\),
  \[\frac{\widetilde B_{\lambda_\kappa}(\kappa,N)}{\widetilde B_0(\kappa,N)}\longrightarrow0 \quad\Longleftrightarrow\quad |\kappa-\lambda_\kappa|=o\!\left(\kappa^{d_*/(n+d_*)}\right).\]
  Moreover, if \(\delta_\kappa=\kappa-\lambda_\kappa=o(\kappa)\), then
  \(\widetilde B_{\lambda_\kappa}/\widetilde B_0\asymp_a\kappa^{-d_*}\br{\delta_\kappa}^{n+d_*}\).
  In particular, bounded mismatch retains the order \(\kappa^{-d_*}\).
  These statements concern the specified upper bounds \(\widetilde B_\lambda\).
\end{corollary}

\begin{proof}
  Since \(\delta\mapsto C_\delta a\) is continuous in \(\Sclass(\R)\), by continuity of the STFT
  and weighted integration we have continuity of \(\delta\mapsto\|C_\delta a\|_{\mathcal M_{s_*}}\).
  Since \(C_\delta a\ne0\), by \eqref{eq:asym} and compactness of bounded parameter intervals we obtain
  \(\|C_\delta a\|_{\mathcal M_{s_*}}\asymp_a\br{\delta}^{n+d_*}\) for all \(\delta\in\R\).

  By the recurrence for \(P_{n,\lambda}\), we have
  \(P_{n,\lambda}(x)=(2\pi i\lambda x)^n+O_{\Omega,n}(|\lambda|^{n-1})\),
  uniformly on \(\Omega\) for \(|\lambda|\ge1\).
  Since \(\|x^n\|_{L^r(\Omega)}>0\), by testing on \(h\equiv1\) and using the triangle inequality we obtain
  \(\|C_\lambda\|_{W^{n,r}\to W^{n,r}}\ge c_{\Omega,n,r}|\lambda|^n\) for sufficiently large \(|\lambda|\).
  Moreover, since \(|C_\lambda1|=1\), we have \(\|C_\lambda\|_{W^{n,r}\to W^{n,r}}\ge1\).
  Thus, by the previously proved upper estimate, we have
  \(\|C_\lambda\|_{W^{n,r}\to W^{n,r}}\asymp\br{\lambda}^n\) for every \(\lambda\in\R\).
  By substituting these comparisons into the definition of \(\widetilde B_\lambda\), we obtain the first assertion.

  Since \(\br{\lambda_\kappa}^n\ge1\), if \(\widetilde B_{\lambda_\kappa}/\widetilde B_0\to0\),
  we have \(\br{\delta_\kappa}/\br{\kappa}\to0\).
  Thus we have \(\delta_\kappa=o(\kappa)\) and \(\lambda_\kappa/\kappa\to1\).
  By the first assertion, whenever \(\delta_\kappa=o(\kappa)\), we obtain
  \(\widetilde B_{\lambda_\kappa}/\widetilde B_0\asymp_a\kappa^{-d_*}\br{\delta_\kappa}^{n+d_*}\).
  Since \(d_*/(n+d_*)<1\), we obtain both directions of the stated equivalence.
  Moreover, by taking \(\delta_\kappa\) bounded in this comparison, we obtain the final assertion.
\end{proof}

%


\section{Numerical Experiments}\label{chap:numerical-experiments}


The theoretical analysis addresses three main questions. First, can an STSAFT-adapted shallow network efficiently approximate signals with pronounced oscillations and quadratic-phase structure? Second, can the approximation behavior predicted by the theoretical $N^{-1/2}$ upper bound, together with the improvement in the bound caused by phase matching, be reflected in computable finite-dimensional examples? Third, does restricting the range of the spatial translation parameter affect global approximation? To investigate these questions from complementary perspectives, we conduct the following four experiments.

Experiment 1 compares a phase-adapted dictionary with a shallow ReLU spline under approximately matched nominal full-network parameter counts; the fitted output dimensions are not matched. Both function-value error and derivative error are measured, with particular attention to how the two dictionaries behave as the quadratic-phase parameter increases. Experiment 2 varies the number of network blocks $N$ and plots the error against $N$ on log-log axes, using a line with slope $-1/2$ as a reference. Experiment 3 varies the quadratic-phase parameter $\lambda$ of the dictionary and compares the cases $\lambda=0$, $\lambda=\kappa/2$, and $\lambda=\kappa$. It also plots a separate proxy for the theoretical upper-bound ratio as a function of phase mismatch. Experiment 4 restricts the translation parameter when using a compactly supported window and employs a smooth, compactly supported target located away from the origin to illustrate the global error floor caused by an insufficient translation range.

The experiments below are numerical diagnostics of the mechanisms suggested by the theory; they do not constitute proofs of the theoretical results. In particular, the $N^{-1/2}$ rate in the paper is an upper bound for the best $N$-term approximation error. It does not imply that the error produced by every training algorithm must decay with exactly this slope. 

\subsection{Dictionaries, Error Metrics, and Numerical Implementation}\label{sec:exp-method}

\subsubsection{A Phase-Adapted Local Dictionary}

Let the Gaussian window be $g(u)=e^{-u^2/2}$ and choose the gated function
\[
 H(t)=t_+e^{-t^2/2},\qquad t_+=\max\{t,0\},
\]
which corresponds to the real-valued network construction considered in the paper. For parameters $y_j,\eta_j,\beta_j\in\mathbb R$, define the local envelope
\[
 h_j(x)=H(\eta_jx+\beta_j)g(x-y_j).
\]
For a fixed-scale quadratic-phase family, let $\lambda\in\mathbb R$ denote the phase parameter of the dictionary and consider the real-valued approximation
\begin{equation}
 v_{N,\lambda}(x)=\sum_{j=1}^{N}h_j(x)
 \bigl[a_j\cos(\pi\lambda x^2)+b_j\sin(\pi\lambda x^2)\bigr].
 \label{eq:exp-real-network}\tag{4.1}
\end{equation}
This representation separates the local envelope from the quadratic phase: the functions $h_j$ model local amplitude and lower-order oscillatory structure, while the sine and cosine factors explicitly encode the quadratic phase. The family in \eqref{eq:exp-real-network} is a fixed-scale subfamily of the real-valued network dictionary discussed in the paper. In the implementation, the hidden parameters are sampled reproducibly and the output coefficients $a_j,b_j$ are then fitted by least squares.

For Experiment 1, we retain the nominal full-network parameter comparison as one diagnostic:
one STSAFT block is counted as five real parameters and a width-$M$ ReLU network
as approximately three per hidden unit, so the nominal comparison uses
$M=\operatorname{round}(5N/3)$. This is not a match of the fitted least-squares
dimensions, because all hidden features are fixed before fitting the output
coefficients. The actual ReLU spline basis used in the code is
\[
  1,\quad (x-\tau_1)_+,\ldots,(x-\tau_M)_+,
\]
where $\tau_1=-2$ and the remaining knots are the first $M-1$ points of
$\operatorname{linspace}(-1,1,M)$ after omitting the left endpoint. Thus the
ReLU design has $M+1$ columns, including the output bias. The STSAFT-inspired
real-valued design has $2N$ columns, namely the cosine and sine branches.
When $\lambda=0$, all sine columns vanish identically, leaving at most $N$
active columns. Hence the nominal parameter comparison does not match the
actual output dimension, and the zero-phase case is rank-deficient by
construction.

We therefore add a second comparison matched by the numerical rank of the
value-and-derivative design matrix. For each design, we stack its sampled
function and derivative matrices, normalize every column to unit Euclidean
norm, and count singular values satisfying
\[
  \sigma_j>10^{-10}\sigma_1.
\]
For each pair $(N,\kappa)$, the number of ReLU knots is then selected so that
the ReLU design has the same numerical rank as the STSAFT-inspired design.
The numerical rank is computed before ridge regularization. This rank-matched
comparison controls the effective output dimension more directly than the
nominal full-network parameter count; it does not match hidden-parameter
optimization, since the features remain fixed in both models.

The random features in Experiment 1 are generated with seed $20261026$:
$y_j$ are sampled uniformly from $[-1.15,1.15]$, $\eta_j$ uniformly from
$[-7,7]$, and values with $|\eta_j|<0.55$ are replaced by $0.55$ or $-0.55$
according to their sign (zero is assigned $0.55$). The offsets are
\[
  \beta_j=1-\eta_jy_j+\varepsilon_j,\qquad
  \varepsilon_j\sim\operatorname{Unif}[-0.65,0.65].
\]
The same nested feature bank is reused for all $\kappa$. The evaluation grid
contains $1601$ equally spaced points on $[-1,1]$. Least squares is solved
after column normalization by augmenting the normalized design with
$\sqrt{10^{-8}}I$ and calling a least-squares solver with
relative cutoff $10^{-10}$. Thus the reported fits use ridge parameter
$10^{-8}$; the numerical-rank diagnostic above uses the unregularized,
column-normalized design and the same relative singular-value threshold.

\subsubsection{Error Metrics}

For a target function $f$ and an approximation $v$, we use the relative $L^2$ error
\begin{equation}
 \mathcal E_{L^2}(f,v)=\frac{\|f-v\|_{L^2(I)}}{\|f\|_{L^2(I)}},
 \label{eq:exp-rel-l2}\tag{4.2}
\end{equation}
and the relative $W^{1,2}$ error
\begin{equation}
 \mathcal E_{W^{1,2}}(f,v)=
 \frac{\bigl(\|f-v\|_{L^2(I)}^2+\|f'-v'\|_{L^2(I)}^2\bigr)^{1/2}}
 {\bigl(\|f\|_{L^2(I)}^2+\|f'\|_{L^2(I)}^2\bigr)^{1/2}}.
 \label{eq:exp-rel-w12}\tag{4.3}
\end{equation}
Unless stated otherwise, the interval in \eqref{eq:exp-rel-l2}--\eqref{eq:exp-rel-w12} is $I=[-1,1]$. Numerical integrals are approximated by the composite trapezoidal rule on a uniform grid. In Experiment 1, output coefficients are fitted separately for the $L^2$ and $W^{1,2}$ objectives, avoiding an opaque trade-off between function-value error and derivative error in a single combined objective. Experiments 2 and 3 use a discrete $W^{1,2}$ least-squares objective and additionally report the corresponding $L^2$ error. Experiment 4 uses the interval $[-4,4]$, which contains the target support and all accessible atom supports for the translation radii considered; its numerical integrals therefore approximate the corresponding global norms. Random features are generated with fixed seeds. The feature-sampling laws, parameter ranges, evaluation-grid sizes, solver, regularization, and numerical-rank tolerance are specified above and in the accompanying Python scripts. The fitting objective uses unweighted stacked sample values (and derivative values for the $W^{1,2}$ fit), whereas the reported error norms are evaluated by the composite trapezoidal rule. Because Experiments 1--3 use finite-interval errors, they provide diagnostics of the approximation mechanisms rather than direct numerical verification of the global Sobolev estimate on $\mathbb R$.

For the representative oscillatory target considered in the paper, we use
\begin{equation}
 f_{\omega,\kappa}(x)=e^{-x^2}\sin(\omega x+\pi\kappa x^2),
 \label{eq:exp-target}\tag{4.4}
\end{equation}
with derivative
\begin{equation}
 f'_{\omega,\kappa}(x)=e^{-x^2}\Bigl[-2x\sin(\theta(x))+
 (\omega+2\pi\kappa x)\cos(\theta(x))\Bigr],\qquad
 \theta(x)=\omega x+\pi\kappa x^2.
 \label{eq:exp-target-derivative}\tag{4.5}
\end{equation}
Using the explicit formula \eqref{eq:exp-target-derivative} avoids the additional error that would be introduced by approximating the target derivative with finite differences.

\subsection{Experiment 1: Phase-Adapted Dictionary versus Shallow ReLU}\label{sec:exp-architecture}


Experiment 1 compares the dictionary in \eqref{eq:exp-real-network} with a shallow ReLU spline for the target in \eqref{eq:exp-target}, measuring both function-value approximation and first-derivative reconstruction.

We use $I=[-1,1]$, fix $\omega=10$, and consider
\[
 \kappa\in\{0,1,2,4,8\},\qquad N\in\{24,48,72\},\qquad
 M=\operatorname{round}(5N/3).
\]
For the phase-adapted dictionary, we set $\lambda=\kappa$. The shallow ReLU baseline uses a linear-spline dictionary with uniformly spaced fixed knots, and its output coefficients are found by least squares. We separately minimize the discrete $L^2$ and $W^{1,2}$ objectives and evaluate the results using \eqref{eq:exp-rel-l2}--\eqref{eq:exp-rel-w12}.

\begin{figure}[htbp]
 \centering
 \includegraphics[width=0.92\linewidth]{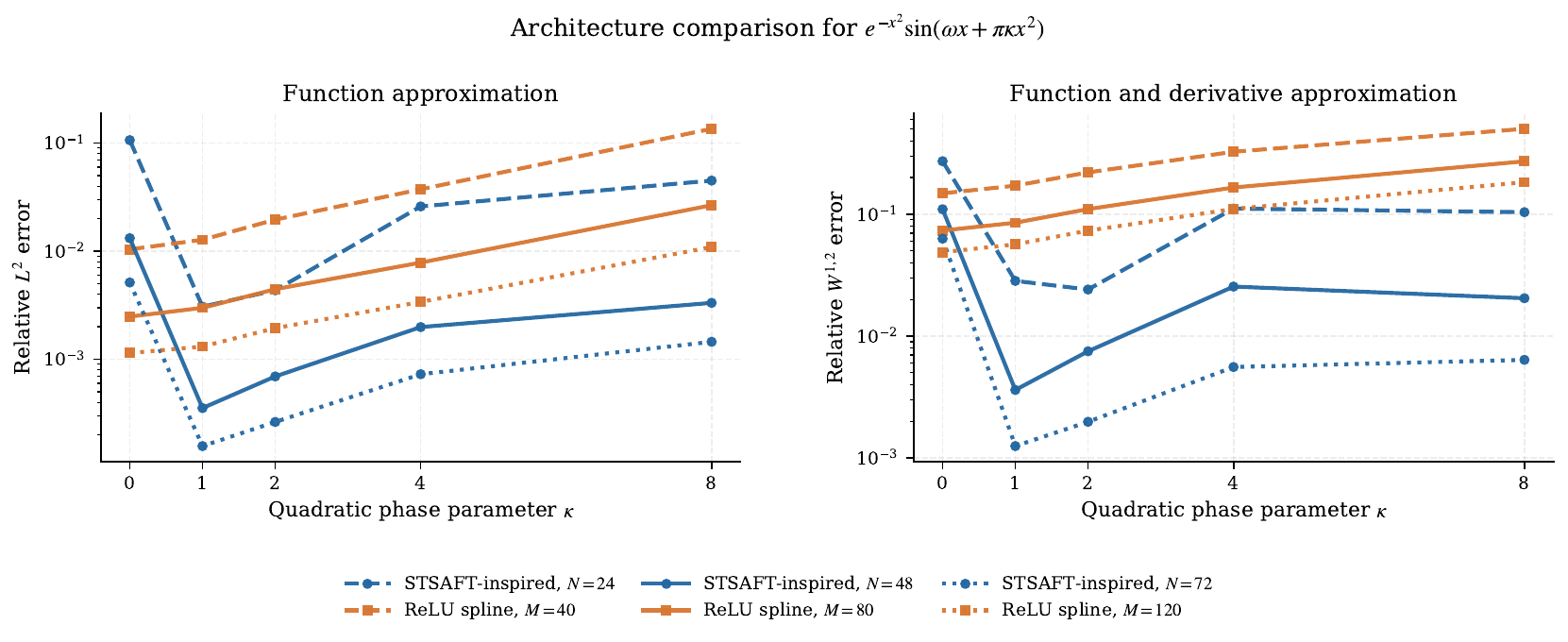}
 \caption{Experiment 1: Relative $L^2$ and $W^{1,2}$ errors for the phase-adapted dictionary and shallow ReLU spline under approximately matched nominal full-network parameter counts; the fitted output dimensions are not matched.}
 \label{fig:exp1-comparison}
\end{figure}

Figure~\ref{fig:exp1-comparison} shows that the errors of both dictionaries generally decrease as the dictionary size increases. For $N=72$, the relative $W^{1,2}$ errors of the STSAFT-inspired dictionary are approximately $1.3\times10^{-3}$, $2.0\times10^{-3}$, $5.6\times10^{-3}$, and $6.4\times10^{-3}$ for $\kappa=1,2,4,8$, respectively. Under the same nominal parameter budget, the corresponding errors of the ReLU spline are approximately $5.7\times10^{-2}$, $7.3\times10^{-2}$, $1.1\times10^{-1}$, and $1.8\times10^{-1}$. In this experiment, the phase-adapted dictionary gives a clear advantage in derivative approximation and retains relatively small errors for stronger quadratic phases.

\begin{figure}[htbp]
 \centering
 \includegraphics[width=0.92\linewidth]{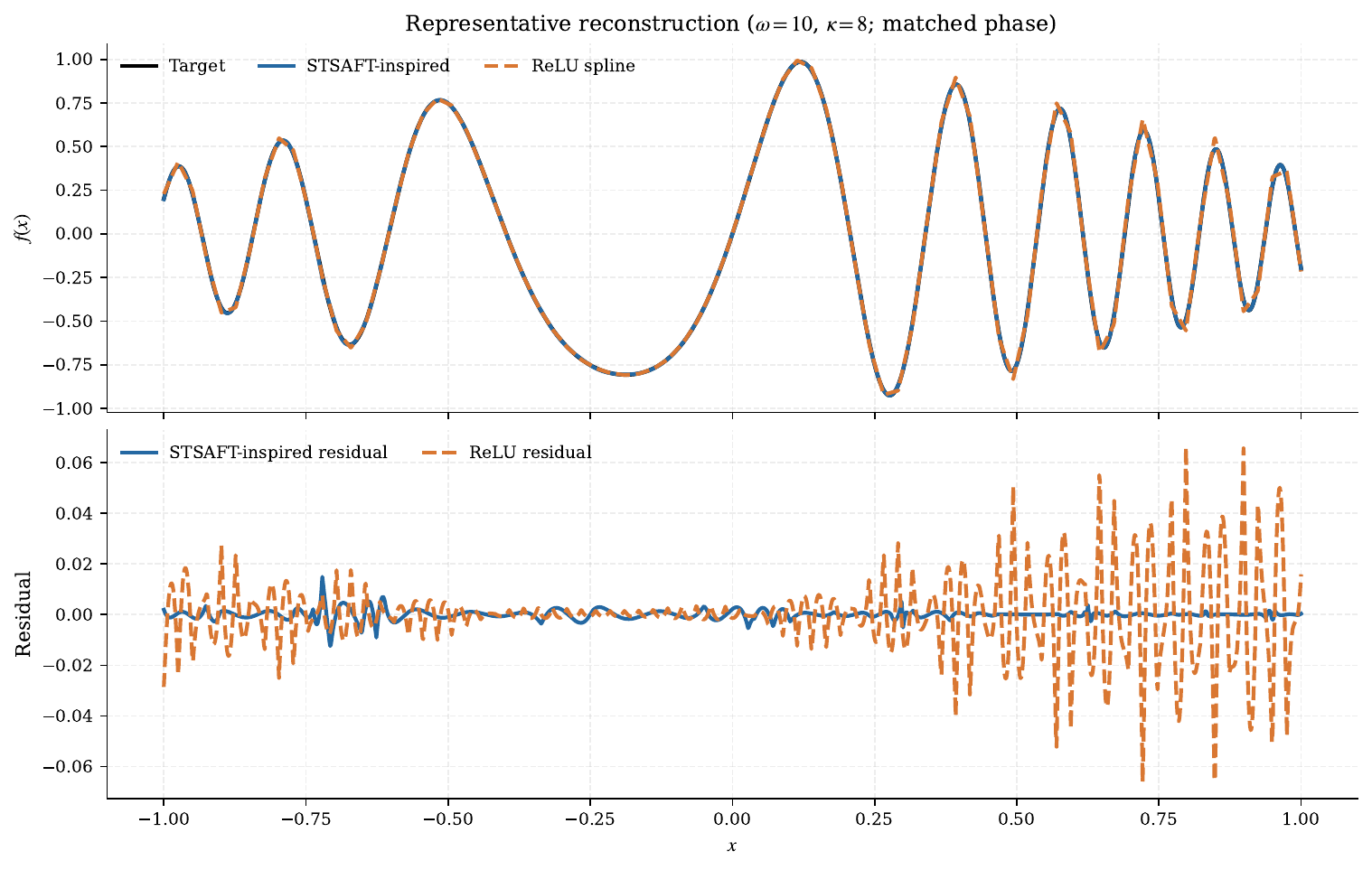}
 \caption{Representative reconstruction and residual in Experiment 1, with $\omega=10$, $\kappa=8$, and $N=48$. The output coefficients in this figure are fitted using the $L^2$ objective to make the function-value fit easier to inspect.}
 \label{fig:exp1-reconstruction}
\end{figure}

For the function-value error at $N=72$, the relative $L^2$ errors of the phase-adapted dictionary for $\kappa=1,2,4,8$ are approximately $1.6\times10^{-4}$, $2.6\times10^{-4}$, $7.3\times10^{-4}$, and $1.5\times10^{-3}$, respectively; all are lower than the corresponding ReLU spline errors. In the nominal comparison, the ReLU spline has smaller errors at $\kappa=0$.
However, this case has only $N$ active STSAFT columns because the sine branch
vanishes. The rank-matched results provide a more informative comparison.
For $N=72$, the numerical ranks of the STSAFT design are $72$ at $\kappa=0$,
$141$ at $\kappa=1$, $142$ at $\kappa=2$, and $144$ at $\kappa=4,8$.
The ReLU widths in the rank-matched control are chosen to give exactly these
same numerical ranks. The complete rank data for Experiment 1 are summarized
in Table~\ref{tab:exp1-design-ranks}; the ReLU design is matched to the
corresponding entry in every row.
\begin{table}[htbp]
 \centering
 \begin{tabular}{c|ccccc}
  $N$ & $\kappa=0$ & $\kappa=1$ & $\kappa=2$ & $\kappa=4$ & $\kappa=8$\\
  \hline
  $24$ & $24$ & $48$ & $48$ & $48$ & $48$\\
  $48$ & $48$ & $96$ & $96$ & $96$ & $96$\\
  $72$ & $72$ & $141$ & $142$ & $144$ & $144$
 \end{tabular}
 \caption{Numerical ranks of the STSAFT-inspired value-and-derivative design
 matrices in Experiment 1. The ReLU spline width is chosen so that its
 numerical rank equals the corresponding entry.}
 \label{tab:exp1-design-ranks}
\end{table}
At nonzero chirp values, the STSAFT-inspired dictionary
still has substantially smaller $W^{1,2}$ errors: for $N=72$, these are
approximately $1.25\times10^{-3}$, $1.98\times10^{-3}$,
$5.58\times10^{-3}$, and $6.39\times10^{-3}$ for
$\kappa=1,2,4,8$, whereas the rank-matched ReLU errors are approximately
$4.81\times10^{-2}$, $6.17\times10^{-2}$,
$9.18\times10^{-2}$, and $1.52\times10^{-1}$.
At $\kappa=0$, the rank-matched $L^2$ error remains smaller for ReLU, while
the $W^{1,2}$ comparison is closer and reverses for $N=48,72$ in this
particular fixed-feature experiment. These observations support a
phase-dependent advantage in the tested nonzero-chirp cases, but do not
establish universal superiority of either dictionary. A representative
reconstruction is shown in Figure~\ref{fig:exp1-reconstruction}.

\begin{figure}[htbp]
 \centering
 \includegraphics[width=0.98\linewidth]{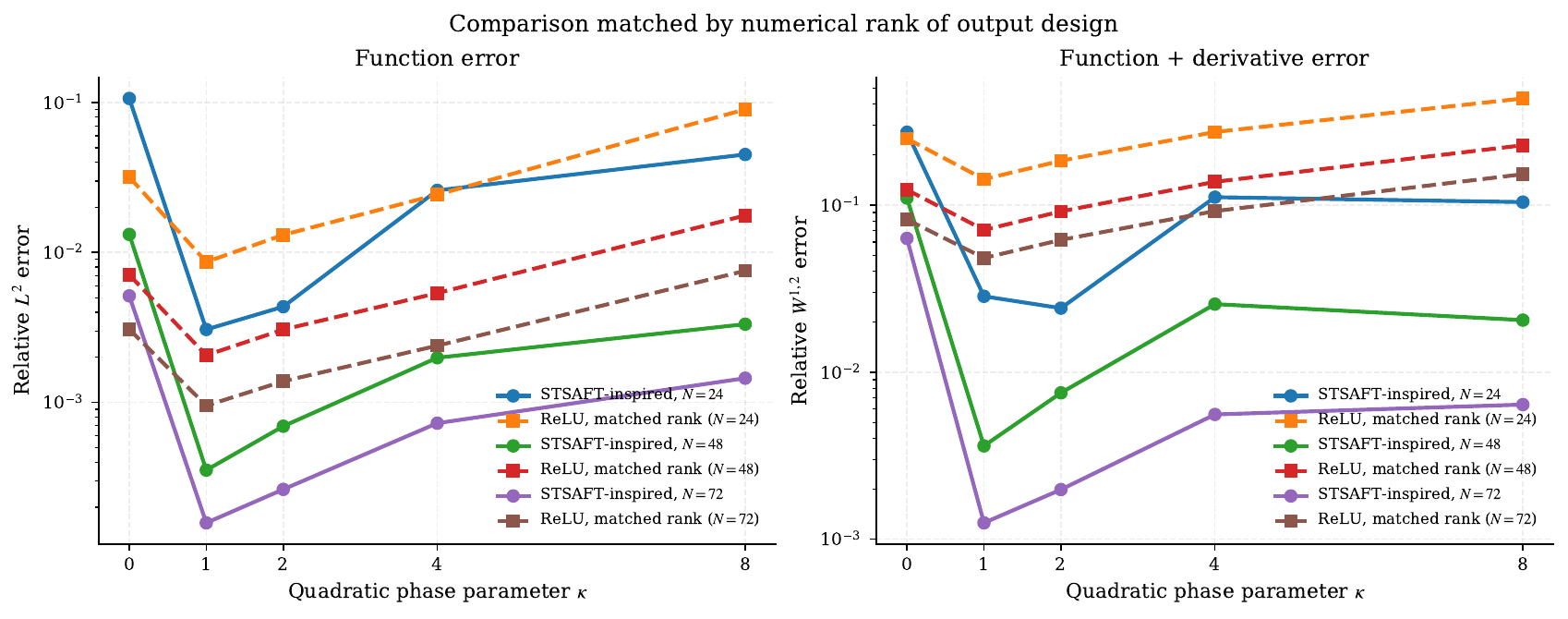}
 \caption{Experiment 1: Comparison after matching the numerical ranks of the
 value-and-derivative design matrices. The numerical rank is computed after
 column normalization with relative singular-value tolerance $10^{-10}$.}
 \label{fig:exp1-effective-dof}
\end{figure}

\subsection[Experiment 2: Empirical Convergence of the N-Term Error]{Experiment 2: Empirical Convergence of the $N$-Term Error}\label{sec:exp-rate}


The global approximation theorem gives an upper bound of the form
\begin{equation}
 \inf_{v_N\in\Sigma_N(\mathcal D_A)}\|f-v_N\|_{W^{1,2}(\mathbb R)}
 \le C N^{-1/2}J_{A,1,2}(f).
 \label{eq:theoretical-rate}\tag{4.6}
\end{equation}
The bound in \eqref{eq:theoretical-rate} is stated on $\mathbb R$, whereas the numerical errors in this experiment are evaluated on $I=[-1,1]$. Thus, the computation is a finite-domain diagnostic of error decay as the dictionary grows; it is not a direct numerical verification of the global estimate. We set $\omega=10$, $\kappa=4$, and consider
\[
 N\in\{8,16,32,64,128,256\}.
\]
For each $N$, we use a nested feature library, so that a dictionary of smaller size consists of the first $N$ atoms of a larger library. The experiment is repeated with six random seeds, $20261026,\ldots,20261031$. For each seed, $y_j$ is uniform on $[-1.18,1.18]$, $\eta_j$ is uniform on $[-7,7]$ with $|\eta_j|<0.55$ clipped to magnitude $0.55$, and $\beta_j=1-\eta_jy_j+\varepsilon_j$ with $\varepsilon_j$ uniform on $[-0.65,0.65]$. The evaluation grid has $1401$ equally spaced points on $[-1,1]$. The output coefficients are fitted by the same column-normalized ridge least-squares procedure with ridge $10^{-8}$ and solver cutoff $10^{-10}$. We also record the numerical rank of each stacked value-and-derivative design, using the column-normalized singular-value threshold $10^{-10}\sigma_1$. The observed ranks range from $16$ for $N=8$ to $384$--$418$ for $N=256$ across the six seeds, showing that the nominal $2N$ columns need not all be numerically independent at the stated tolerance. We plot the median error and interquartile range, and include a reference line with slope $-1/2$ in the log-log plot of the $W^{1,2}$ error. We also fit an empirical slope using
\[
 \log E_N=\alpha+\beta\log N+\varepsilon_N.
\]

\begin{figure}[htbp]
 \centering
 \includegraphics[width=0.98\linewidth]{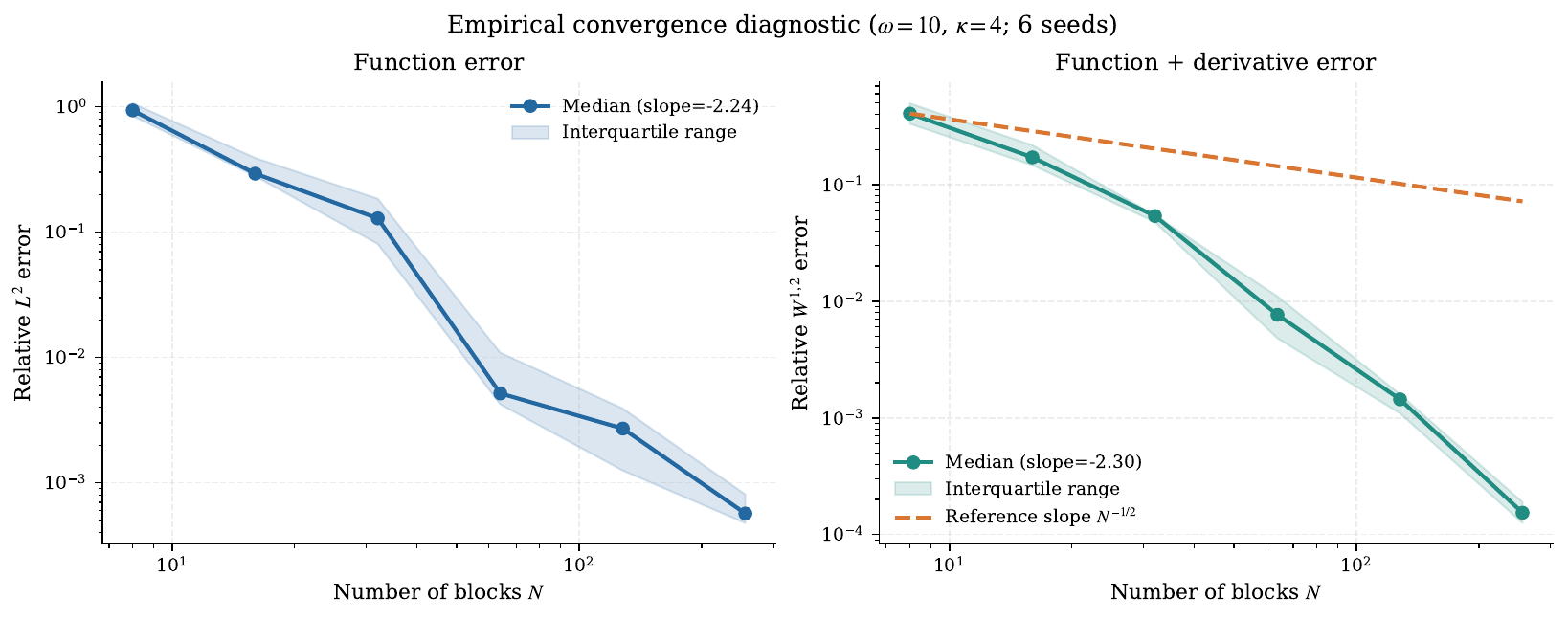}
 \caption{Experiment 2: Median relative errors, interquartile ranges over six random seeds, and the $N^{-1/2}$ reference slope.}
 \label{fig:exp2-rate}
\end{figure}

Figure~\ref{fig:exp2-rate} presents the log-log curves for the function error and the Sobolev error. For the current seeds and discretization, the fitted slopes over the tested range are approximately $-2.24$ for the relative $L^2$ error and $-2.30$ for the relative $W^{1,2}$ error. In this example, the empirical errors decrease faster than the $N^{-1/2}$ reference line. This does not contradict \eqref{eq:theoretical-rate}, because that estimate is a general upper bound and does not assert that the actual error must decay at exactly the rate $N^{-1/2}$. Moreover, the present experiment uses a fixed random dictionary and only a finite set of values of $N$, so the observed slopes do not prove the theoretical best-approximation rate. A more cautious conclusion is that, for this smooth target with a phase-adapted structure, the numerical error decreases rapidly as the dictionary grows.

\subsection{Experiment 3: Quadratic-Phase Matching and Mismatch Scale}\label{sec:exp-phase}


The phase-matching result concerns the relation between the phase of the target and the phase of the dictionary. Define
\[
 C_\lambda u(x)=e^{i\pi\lambda x^2}u(x),\qquad
 f_\kappa(x)=C_\kappa a(x),\qquad a(x)=e^{-x^2}.
\]
We use a fixed collection of local envelope atoms to construct $\mathcal D_\lambda=C_\lambda\mathcal D_0$ and consider
\[
 \kappa\in\{1,2,4,8,12,16\},\qquad
 \lambda\in\{0,\kappa/2,\kappa\},\qquad N=64.
\]
To track the effect of $C_\lambda$ directly, we first perform the computation with a complex-valued dictionary. Each of the $N$ complex output coefficients represents two real scalar degrees of freedom, so the fitted output space has at most $2N$ real dimensions. For the current samples, the complex numerical rank is $64$ for each tested $(\kappa,\lambda)$, corresponding to $128$ effective real output degrees of freedom. The rank is computed from the column-normalized complex stacked value-and-derivative design, with relative singular-value tolerance $10^{-10}$. The real-valued rank-loss observation for the sine/cosine form in Experiment 1 does not apply to this complex formulation. The random-feature seed is $20261026$, the spatial centers are uniform on $[-1.15,1.15]$, the slopes are uniform on $[-7,7]$ with magnitude clipped below $0.55$, and the offsets follow $\beta_j=1-\eta_jy_j+\varepsilon_j$, $\varepsilon_j\sim\operatorname{Unif}[-0.65,0.65]$. The evaluation grid has $1601$ equally spaced points on $[-1,1]$. The complex least-squares solver uses column normalization, ridge $10^{-8}$, and relative cutoff $10^{-10}$. We record both relative $L^2$ and relative $W^{1,2}$ errors.

In addition, motivated by the structure of the theoretical upper bound, define the dimensionless proxy
\begin{equation}
 Q_\kappa^{\mathrm{pow}}(\lambda)=
 \frac{\langle\lambda\rangle^n\langle\kappa-\lambda\rangle^{n+3/2-1/r}}
 {\langle\kappa\rangle^{n+3/2-1/r}},
 \qquad \langle t\rangle=(1+t^2)^{1/2},\quad n=1,\quad r=2.
 \label{eq:bound-ratio-proxy}\tag{4.7}
\end{equation}
We set $\lambda_\kappa=\kappa-\kappa^\alpha$ and consider $\alpha=0.2,0.4,0.5,0.6,0.8$. The resulting values of $Q_\kappa^{\mathrm{pow}}(\lambda_\kappa)$ are examined as functions of $\kappa$. Note that \eqref{eq:bound-ratio-proxy} is a proxy for the ratio of phase-dependent factors in the upper bound. It is not the actual approximation error, nor is it an exact ratio that includes all constants.

\begin{figure}[htbp]
 \centering
 \includegraphics[width=0.98\linewidth]{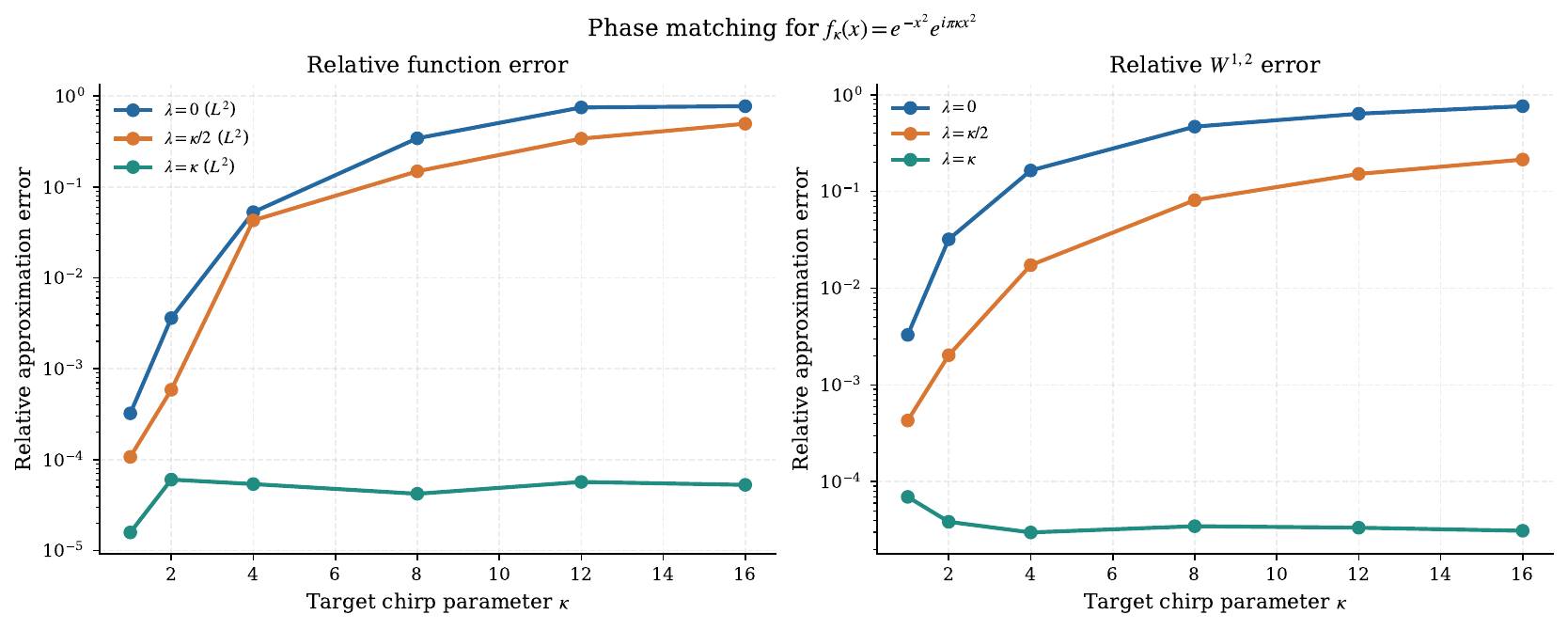}
 \caption{Experiment 3: Relative $L^2$ and $W^{1,2}$ errors for different dictionary phase parameters $\lambda$.}
 \label{fig:exp3-phase}
\end{figure}

Figure~\ref{fig:exp3-phase} shows that when $\lambda=\kappa$, the error is relatively insensitive to an increase in $\kappa$. In this experiment, for $\kappa\geq 4$, the relative $L^2$ error in the phase-matched case is approximately between $4\times10^{-5}$ and $6\times10^{-5}$. In contrast, the unmatched case $\lambda=0$ deteriorates substantially as $\kappa$ increases: at $\kappa=16$, its relative $L^2$ error is approximately $0.77$ and its relative $W^{1,2}$ error is approximately $0.76$. The intermediate phase $\lambda=\kappa/2$ typically performs between these two cases. This behavior is consistent with the theoretical motivation that phase matching simplifies the envelope to be approximated, but the conclusion is specific to the dictionary and targets used here.

\begin{figure}[htbp]
 \centering
 \includegraphics[width=0.88\linewidth]{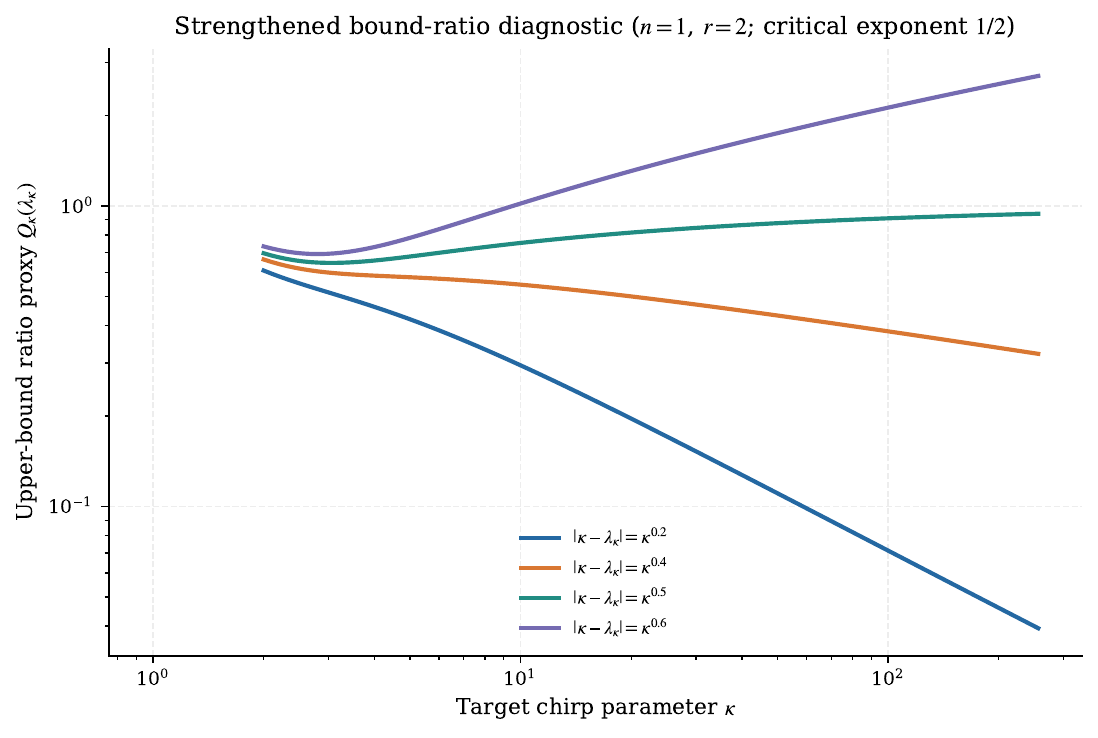}
 \caption{Experiment 3: For $n=1$, the proxy for the theoretical upper-bound ratio as a function of $\kappa$, when the phase mismatch is $|\kappa-\lambda_\kappa|=\kappa^\alpha$. The critical exponent is $1/2$ for $n=1,r=2$.}
 \label{fig:exp3-mismatch}
\end{figure}

For the power proxy in \eqref{eq:bound-ratio-proxy}, when $\lambda_\kappa=\kappa-\kappa^\alpha$ and $n=1$, its asymptotic power under the strengthened estimate is $Q_\kappa^{\mathrm{pow}}\asymp\kappa^{2\alpha-1}$. Thus, this proxy tends to zero if $\alpha<1/2$, remains of constant order if $\alpha=1/2$, and grows if $\alpha>1/2$. This exponent is the mismatch scale of the strengthened specified upper bound; it is not a transition in the empirical errors. Figure~\ref{fig:exp3-mismatch} illustrates this upper-bound scale diagnostic. The exponent $1/2$ pertains to the strengthened specified upper-bound ratio, not to a necessary transition point in empirical training errors.

\subsection{Experiment 4: Global Approximation Obstruction under Restricted Translations}\label{sec:exp-translation}


Suppose that the window $g$ has compact support and the dictionary only permits translation parameters $y\in[-R,R]$. Then every atom $g(x-y)$ is supported in $[-R-r_g,R+r_g]$, where $\operatorname{supp}g\subset[-r_g,r_g]$. Denote this reachable region by
\[
 K_R=[-R-r_g,R+r_g].
\]
Any finite network composed of such atoms vanishes on $K_R^c$. Consequently, every target $f$ satisfies the unavoidable lower bound
\begin{equation}
 \|f-v_N\|_{L^2(\mathbb R)}\geq
 \|f\,\mathbf 1_{K_R^c}\|_{L^2(\mathbb R)},
 \qquad v_N\in\Sigma_N(\mathcal D_R).
 \label{eq:translation-lower-bound}\tag{4.8}
\end{equation}

To visualize this obstruction, we use a smooth compactly supported bump function whose support is located away from the origin, with center $2.65$ and support radius $0.38$. The compactly supported window has radius $r_g=0.5$. We choose
\[
 R\in\{0.5,1,1.5,2,2.5,3,3.5\}.
\]
For each value of $R$, we use $N=180$ local atoms. Their centers are equally spaced on $[-R,R]$; the slopes are sampled uniformly on $[-6,6]$ and clipped away from zero to magnitude at least $0.6$, and the offsets are $\beta_j=1-\eta_jy_j+\varepsilon_j$ with $\varepsilon_j$ uniform on $[-0.65,0.65]$. The random seed is $20261026$. The target is the bump $f(x)=\exp[-1/(1-s^2)]$ for $|s|<1$ and zero otherwise, where $s=(x-2.65)/0.38$. The window is the same bump profile centered at zero with radius $0.5$. The grid has $2401$ equally spaced points on $[-4,4]$. The function and derivative are fitted by column-normalized ridge least squares with ridge $10^{-8}$ and relative solver cutoff $10^{-10}$. Each of the seven stacked designs has numerical rank $180$ at relative rank tolerance $10^{-10}$. We compare the observed relative error with the relative $L^2$ lower bound in \eqref{eq:translation-lower-bound}.

\begin{figure}[htbp]
 \centering
 \includegraphics[width=0.98\linewidth]{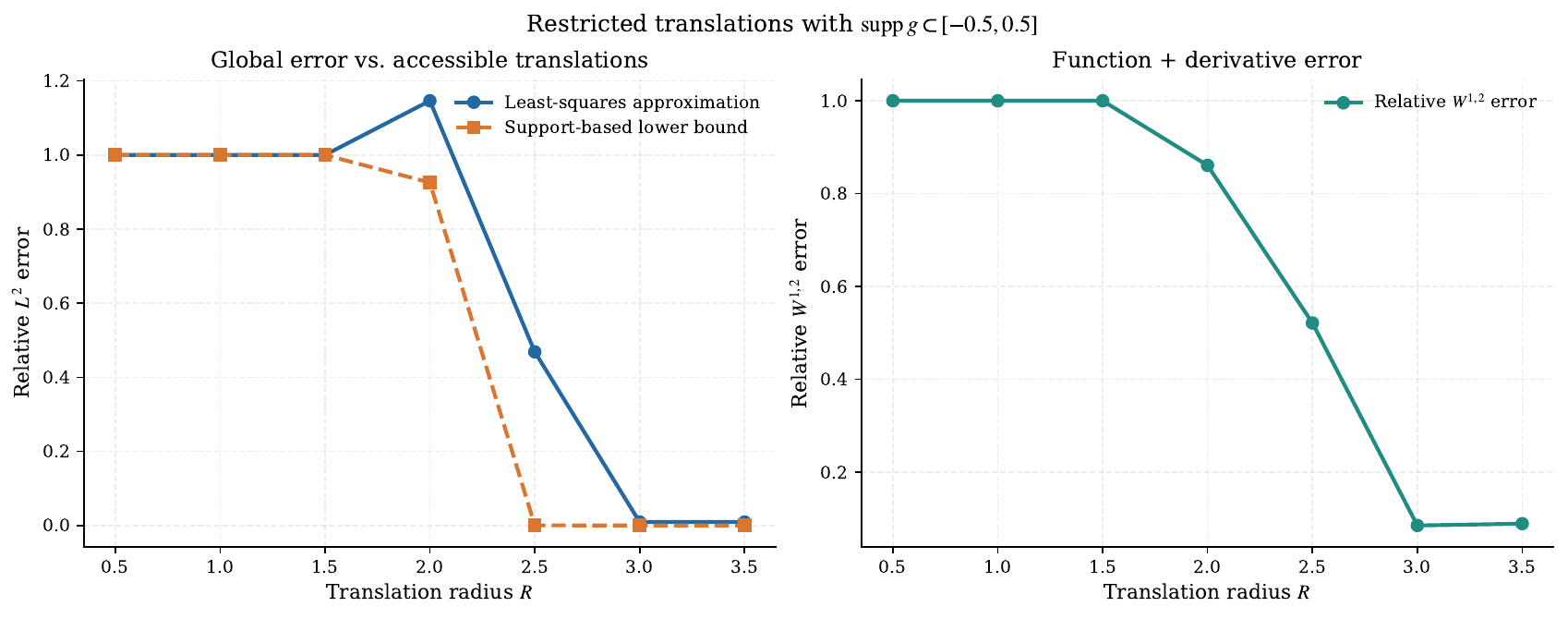}
 \caption{Experiment 4: Global relative $L^2$ error, support-based lower bound, and relative $W^{1,2}$ error as the allowed translation radius varies.}
 \label{fig:exp4-range}
\end{figure}
\begin{figure}[htbp]
 \centering
 \includegraphics[width=0.94\linewidth]{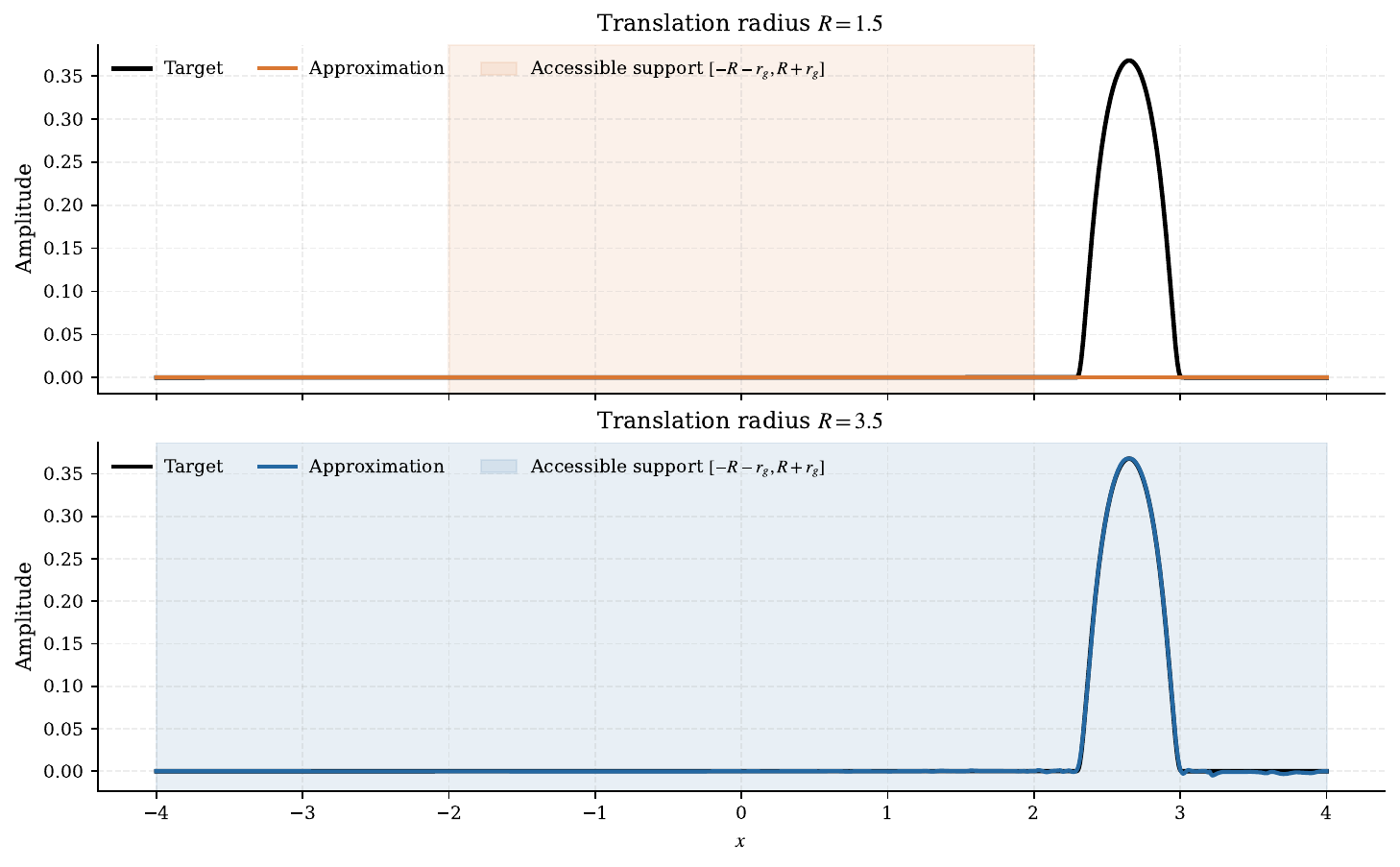}
 \caption{Experiment 4: Target, approximation, and reachable support region for translation radii $R=1.5$ and $R=3.5$.}
 \label{fig:exp4-support}
\end{figure}

Figure~\ref{fig:exp4-range} shows that when $R\leq 1.5$, the reachable region is disjoint from the support of the target, and the network's relative $L^2$ error is exactly $1$, in agreement with \eqref{eq:translation-lower-bound}. At $R=2$, the reachable region covers only part of the target support; the support-based lower bound is approximately $0.926$ and the observed relative $L^2$ error is approximately $1.15$. An error larger than one is possible here because the coefficients are fitted using the $W^{1,2}$ objective rather than by direct minimization of the $L^2$ error. When $R$ increases to $3$, the target support is contained in the reachable region and the observed relative $L^2$ error drops to approximately $9.4\times10^{-3}$. At $R=2.5$, although the support-based lower bound is already close to zero, the finite dictionary still has an empirical error of approximately $0.47$. This illustrates that covering the target support removes the support-based obstruction, but does not by itself guarantee high-accuracy approximation with a finite number of atoms. Figure~\ref{fig:exp4-support} further compares reconstructions for a small and a larger translation range.

\subsection{Summary}\label{sec:exp-summary}

The four experiments investigate the structural factors highlighted by the theory from four perspectives: model comparison, error decay with dictionary size, phase matching, and the range of translation parameters. Experiment 1 indicates that, after matching the numerical output ranks, the STSAFT-inspired dictionary remains substantially more effective for derivative reconstruction in the tested nonzero-chirp cases; the zero-phase case also illustrates why the sine-branch rank loss must be accounted for. Experiment 2 shows that the error decreases as the number of atoms increases, while the observed empirical slope should not be interpreted as a proof of the theoretical rate. Experiment 3 demonstrates a substantial reduction in approximation error for the current complex-valued test when the phases are matched, and separately uses an upper-bound proxy to illustrate the critical mismatch scale. Experiment 4 verifies the global support obstruction caused jointly by a compactly supported window and a restricted translation range.


\bibliographystyle{plainnat}
\bibliography{thesis}
\section*{Data Availability}
Data sharing is not applicable because no datasets were generated or analyzed in this theoretical study.

\end{document}